\documentclass[11pt]{amsart}
\usepackage{epsfig}
\input rlepsf

\newcommand{\bfz}{{\mathbb {Z}}}

\newcommand{\bfc}{{\mathbb {C}}}
\newcommand{\bfn}{{\mathbb {N}}}
\newcommand{\bfr}{{\mathbb {R}}}

\renewcommand{\t}{\mathbf t}
\newcommand{\s}{\mathbf s}
\newcommand{\Li}{\mathbb {L}}
\newcommand{\tb}{\mbox{\rm tb}}
\newcommand{\lk}{{\ell k}}
\newcommand{\rot}{\mbox{\rm rot}}

\begin{document}

\title{Surgery diagrams for contact $3$-manifolds}
\author{Fan Ding}
\author{Hansj\"org Geiges}
\author{Andr\'{a}s I. Stipsicz}
\address{Department of Mathematics, Peking University, Beijing 100871, 
P.R. China}
\email{dingfan@math.pku.edu.cn}
\address{Mathematisches Institut, Universit\"at zu K\"oln, 
Weyertal 86-90, 50931 K\"oln, Germany}
\email{geiges@math.uni-koeln.de}
\address{A. R\'enyi Institute of Mathematics, Hungarian Academy of Sciences,
Budapest, Re\'altanoda utca 13--15, H-1053, Hungary}
\email{stipsicz@math-inst.hu}
\date{}
\begin{abstract}
In two previous papers, the two first-named authors introduced a notion
of contact $r$-surgery along Legendrian knots in contact $3$-manifolds.
They also showed how (at least in principle) to convert any contact
$r$-surgery into a sequence of contact $(\pm 1)$-surgeries, and used this
to prove that any (closed) contact $3$-manifold can be obtained
from the standard contact structure on $S^3$ by a sequence of such
contact $(\pm 1)$-surgeries.

In the present paper, we give a shorter proof of that result and
a more explicit algorithm for turning a contact $r$-surgery into
$(\pm 1)$-surgeries. We use this to give explicit surgery diagrams
for all contact structures on $S^3$ and $S^1\times S^2$, as well
as all overtwisted contact structures on arbitrary closed,
orientable $3$-manifolds. This amounts to a new proof of the Lutz-Martinet
theorem that each homotopy class of $2$-plane fields on such a manifold
is represented by a contact structure.
\end{abstract}
\maketitle

\newtheorem{theorem}{Theorem}[section]
\newtheorem{lemma}[theorem]{Lemma}
\newtheorem{proposition}[theorem]{Proposition}
\newtheorem{corollary}[theorem]{Corollary}
\newtheorem{conjecture}[theorem]{Conjecture}
\newtheorem{question}[theorem]{Question}

\theoremstyle{definition}
\newtheorem{definition}[theorem]{Definition}
\newtheorem{example}[theorem]{Example}
\newtheorem{examples}[theorem]{Examples}
\newtheorem{xca}[theorem]{Exercise}
\newtheorem{remark}[theorem]{Remark}
\newtheorem{remarks}[theorem]{Remarks}

\theoremstyle{remark}

\numberwithin{equation}{section}

\newcommand{\blankbox}[2]{%
  \parbox{\columnwidth}{\centering
    \setlength{\fboxsep}{0pt}%
    \fbox{\raisebox{0pt}[#2]{\hspace{#1}}}%
  }%
}

\section{Introduction}
\label{first}
Let $Y$ be a closed, orientable $3$-manifold. A {\em coorientable contact
structure} on $Y$ is the kernel $\xi =\ker\alpha$ of a differential
$1$-form on $Y$ with the property that $\alpha\wedge d\alpha$ is
a volume form. Fixing a coorientation of $\xi$ amounts to fixing
$\alpha$ up to multiplication with a postive function. In the sequel,
we shall assume implicitly that our contact structures are
cooriented; moreover, we equip $Y$ with the orientation induced by the volume
form $\alpha\wedge d\alpha$. This ensures that when below we realise
certain $(Y,\xi )$ as the boundary of an almost complex $4$-manifold
$(X,J)$, the orientation of $Y$ induced by $\xi$ coincides with
the orientation of $Y$ as the boundary of the manifold $X$ (oriented by~$J$).

The standard contact structure $\xi_{st}$ on the $3$-sphere
$S^3\subset {\mathbb R}^4$ (with cartesian coordinates $x,y,z,t$)
is defined as the kernel of
\[ \alpha_{st}=x\, dy-y\, dx+z\, dt-t\, dz \]
or, equivalently, as the complex tangencies of $S^3\subset{\mathbb C}^2$.
For other basics of contact geometry we refer to
\cite{elia92}; for Legendrian knots and their presentation via
front projections see \cite{gomp98}, \cite{etny}; for the general 
differential topological background of contact geometry
see~\cite{geig}.

A {\em Legendrian knot} $K$ in a contact $3$-manifold $(Y,\xi )$ is
a knot that is everywhere tangent to~$\xi$. Such knots come with
a canonical {\em contact framing}, defined by a vector field along $K$
that is transverse to~$\xi$. Recall that $(Y,\xi )$ is called {\em overtwisted}
if it contains an embedded disc $D^2\subset Y$ with boundary $\partial D^2$
a Legendrian knot whose contact framing equals the framing it
receives from the disc~$D^2$. If no such disc exists, the contact structure
is called {\em tight}.

In~\cite{dige01} a notion of contact
$r$-surgery along a Legendrian knot $K$ in a contact
manifold $(Y,\xi )$ was described: This amounts to a topological surgery,
with surgery coefficient $r\in{\mathbb Q}\cup{\infty}$ measured
relative to the contact framing. A contact structure on the surgered
manifold
\[ (Y-\nu K)\cup (S^1\times D^2),\]
with $\nu K$ denoting a tubular neighbourhood of~$K$,
is defined, for $r\neq 0$, by requiring this contact structure to coincide
with $\xi$ on $Y-\nu K$ and its extension over $S^1\times D^2$ to
be tight (on $S^1\times D^2$, not necessarily the whole surgered manifold).
According to~\cite{hond00}, such an extension always exists and is unique
(up to isotopy) for $r=1/k$ with $k\in\bfz$. (For $r=0$, that extension
is necessarily overtwisted and thus requires a different
treatment. For that reason we shall not discuss the case of $0$-surgery
any further in the present paper.) Therefore,
if $r=1/k$ with $k\in{\mathbb Z}$,
there is a canonical procedure for this surgery, that is, the resulting
contact structure on the surgered manifold
is completely determined by the initial manifold
$(Y,\xi)$, the Legendrian knot $K$ in~$Y$, and the surgery coefficient
$r=1/k$.

A contact $(-1)$-surgery corresponds to a symplectic handlebody
surgery in the sense of \cite{elia90}, \cite{wein91} (cf.\ also
Remark~\ref{rem:pm} below). For future reference we record the
following lemma, see~\cite[Prop.~8]{dige01}, \cite[Section~3]{dige}:

\begin{lemma}
\label{lem:cancel}
Contact $(-1)$-surgery along a Legendrian knot $K\subset (Y,\xi )$
and contact $(+1)$-surgery along a Legendrian push-off of $K$
cancel each other.
\end{lemma}

In \cite{dige} the following has been proved:
\begin{theorem}[\cite{dige}]\label{main}
Every (closed, orientable) contact $3$-manifold $(Y, \xi )$ can be 
obtained via contact $(\pm 1)$-surgery on a Legendrian link in
$(S^3, \xi _{st})$.
\end{theorem}

A simple way of proving this theorem relies on the following 
result of Etnyre and Honda:

\begin{theorem}[\cite{etho01}] \label{eh}
Let $(Y_i, \xi _i)$ ($i=1,2$) be two given contact $3$-manifolds and suppose
that $(Y_1, \xi _1)$ is overtwisted. Then there is a Legendrian link
${\mathbb {L}}\subset (Y_1, \xi _1)$ such that contact $(-1)$-surgery on 
${\mathbb {L}}$ produces $(Y_2, \xi _2)$.
\qed
\end{theorem}

\begin{proof}[Proof of Theorem~\ref{main}] Let $(Y_2,\xi_2)=(Y, \xi )$
be given. Let $(Y_1,\xi_1)$ be the contact manifold obtained by
contact $(+1)$-surgery on the Legendrian knot $K$ in
$(S^3,\xi_{st})$ shown in Figure~\ref{shark}. 

\begin{figure}[h]
\centerline{\relabelbox\small
\epsfxsize 8cm \epsfbox{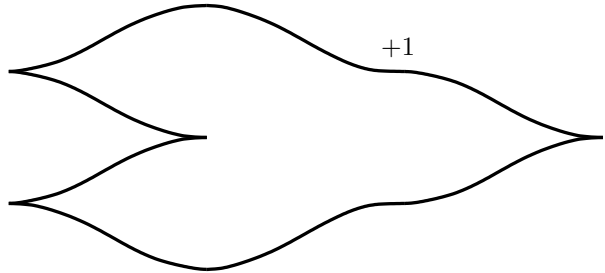}
\extralabel <-3cm,2.9cm> {$+1$}
\endrelabelbox}
  \caption{Legendrian knot producing an overtwisted $S^3$.}
  \label{shark} 
\end{figure} 

That Legendrian knot $K$ has Thurston-Bennequin invariant~$-2$, that is,
the longitude $\lambda_c$ given by the contact framing is related
(homotopically) to
the meridian $\mu$ and standard longitude $\lambda$ of~$K$
(with linking number $\ell k\, (\lambda ,K)=0$) by
$\lambda_c=\lambda -2\mu$. Thus, contact $(+1)$-surgery along $K$ means
that we cut out a tubular neighbourhood of $K$ and glue in a solid
torus by sending its meridian to $\lambda_c+\mu =\lambda-\mu$,
which amounts to a topological $(-1)$-surgery with respect to
the standard framing given by~$\lambda$. Such a surgery is topologically
trivial, that is, $Y_1=S^3$.

\begin{figure}[h]
\centerline{\relabelbox\small
\epsfxsize 12cm \epsfbox{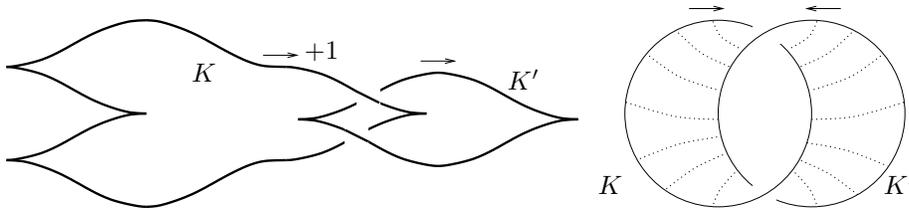}
\extralabel <-8cm,2cm> {$+1$}
\extralabel <-9.5cm,1.7cm> {$K$}
\extralabel <-5.3cm,1.6cm> {$K'$}
\extralabel <-4.1cm,.2cm> {$K$}
\extralabel <-.3cm,.2cm> {$K'$}
\endrelabelbox}
  \caption{The overtwisted disc in $(S^3, \xi _1$).}
  \label{ot} 
\end{figure} 

Figure~\ref{ot} shows that $(S^3, \xi _1)$ is overtwisted:
The surface framing of $K$ determined by the Seifert surface $\Sigma$ of the
Hopf link $K\sqcup K'$ shown in that figure is~$-1$, hence equal to the
framing used for the surgery. This implies that the new meridional disc
$D_m$ in the surgered manifold and $\Sigma$ glued together
define an embedded disc $D_0=D_m\cup_K\Sigma$ in the surgered manifold.
The surface framing of $K'$ determined by $D_0$ is $-1$, which equals the
contact framing of~$K'$. Hence $D_0$ is an overtwisted disc.

It follows that Theorem~\ref{eh} applies and yields the desired 
surgery presentation.
\end{proof}

Notice that, in fact, we have obtained a slightly stronger statement:

\begin{corollary}
Let $(Y, \xi )$ be a contact $3$-manifold. Then there is a Legendrian link 
${\mathbb {L}}\subset (S^3, \xi _{st})$ and a Legendrian knot 
$K\subset (S^3, \xi _{st})$ disjoint from ${\mathbb L}$
such that contact $(+1)$-surgery
on $K$ and contact $(-1)$-surgery on ${\mathbb {L}}$ yield $(Y, \xi )$.
\qed
\end{corollary}

In other words, we can assume that in the surgery presentation  we have a 
single knot on which we do contact $(+1)$-surgery. As the proof shows,
this $K\subset (S^3, \xi _{st})$ can be chosen arbitrarily as long as
$(+1)$-surgery on it results in an overtwisted structure. Needless to say,
different choices for $K$ necessitate different Legendrian links 
${\mathbb {L}}$ for the $(-1)$-surgeries.

\begin{corollary}
For a contact $3$-manifold $(Y, \xi )$ there is a Legendrian knot 
$K^*$ such that $(Y-\nu K^*, \xi \vert_{Y-\nu K^*})$, the complement
of a tubular neighbourhood $\nu K^*$ of~$K^*$, embeds into a 
Stein fillable contact $3$-manifold. In particular,
$(Y-\nu K^*, \xi \vert _{Y-\nu K^*})$ is tight.
\end{corollary}

\begin{proof}
Let $(Y',\xi')$ be the contact manifold obtained by performing
the contact $(-1)$-surgeries along~${\mathbb L}$. This is a Stein fillable
manifold. Our manifold $(Y,\xi )$ is obtained from $(Y',\xi')$ by
a contact $(+1)$-surgery along $K$ (which we may regard as a
Legendrian knot in~$(Y',\xi ')$), that is,
\[ (Y,\xi )=(Y'-\nu K,\xi'\vert_{Y'-\nu K})\cup (S^1\times D^2),\]
where $\xi$ is defined by the unique extension of $\xi'$ over
$S^1\times D^2$ as a tight contact structure on that solid torus.
For a contact $(+1)$-surgery, that contact structure on $S^1\times D^2$
is the unique contact structure on the tubular neighbourhood $\nu K^*$ of
a Legendrian knot~$K^*$. So we may think of $K^*$ as a Legendrian
knot in $(Y,\xi )$ and identify $(Y-\nu K^*, \xi \vert _{Y-\nu K^*})$
with $(Y'-\nu K,\xi'\vert_{Y'-\nu K})$.
\end{proof}

\begin{remark}
The proof of Theorem~\ref{eh} proceeds roughly as follows: If
$(Y_2,\xi _2)$ is also overtwisted, then any 4-dimensional cobordism
from $Y_1$ to $Y_2$ involving only 2-handles can be equipped with a
Stein structure, providing a suitable Legendrian link ${\mathbb L}$
in $(Y_1, \xi_1)$. (Here we use Eliashberg's classification of overwisted
contact structures \cite{elia89} together with his results on the
existence of Stein structures on cobordisms \cite{elia90}.)

For the general case,
consider $(Y_2, \xi _2)\# (S^3, \xi_1)$ (which can be obtained by
performing $(+1)$-surgery on a copy of the knot of Figure~\ref{shark} in a
Darboux chart of $(Y_2, \xi_2)$). Apply the above argument to that
manifold to obtain a Legendrian link ${\mathbb L}'\subset (Y_1,\xi_1)$
such that contact $(-1)$-surgery on ${\mathbb L}'$ yields
$(Y_2,\xi_2)\# (S^3,\xi_1)$.

By Lemma~\ref{lem:cancel}, that first contact $(+1)$-surgery 
can be inverted by a contact $(-1)$-surgery along a suitable 
Legendrian knot~$K^*$, which we may think of as a knot in $(Y_1,\xi_1)$
disjoint from~${\mathbb L}'$. Then ${\mathbb L}={\mathbb L}'\,\sqcup\, K^*$
is the desired link.
\end{remark}

Theorem~\ref{main} can be proved more directly by first reducing it to
the case of overtwisted contact structures on $S^3$, and then giving explicit
surgery diagrams for those structures. Here we shortly describe
this reduction, the explicit diagrams for $S^3$ will be exhibited in
Section~\ref{fourth}.  For the reduction consider, once again,
the manifold $(Y, \xi )\# (S^3,\xi_1)$ constructed via a contact
$(+1)$-surgery on $(Y,\xi )$. It is known that
$Y$ contains a smooth link on which smooth integral surgery provides
$S^3$. Isotoping the components of this link in the {\it overtwisted}
contact 3-manifold $(Y,\xi )\# (S^3, \xi _1)$ we can find, by~\cite{elia90},
a Legendrian link such that contact $(+1)$-surgery on it yields $S^3$
with some contact structure $\xi _r$. (In an overtwisted contact
manifold one can add arbitrary positive or negative twists to the
contact framing of a given Legendrian knot by a suitable band sum with
the boundary of an overtwisted disc.) By taking an additional $(S^3,
\xi _1)$-summand for the whole process, if necessary, we can arrange that
$(S^3, \xi_r)$ is overtwisted.

By inverting the contact $(+1)$-surgeries we
end up with a Legendrian link in $(S^3, \xi _r)$, contact $(-1)$-surgery
on which yields $(Y, \xi)$. This time, however, we do not have
any control on the contact structure $\xi _r$ --- besides it being 
overtwisted. With the help of Eliashberg's classification of overtwisted
contact structures (applied now for $Y=S^3$ only), together with
the mentioned results of Section~\ref{fourth}, we get an alternative proof of
Theorem~\ref{main}.
\subsection*{The algorithm}
In \cite{dige} an algorithm was described (though not entirely
explicitly) for turning a rational contact $r$-surgery into a sequence
of contact $(\pm 1)$-surgeries. Here we extract the relevant information
from~\cite{dige} to formulate an algorithm directly applicable
to a given rational surgery diagram. This algorithm naturally bears
some resemblance to considerations in~\cite{gomp98}. For applications of
this algorithm to the construction of interesting tight contact
structures (e.g.\ ones that are not symplectically semi-fillable)
see~\cite{lista} and~\cite{listb}.

\vspace{1mm}

\begin{figure}[h]
\centerline{\relabelbox\small
\epsfxsize 8cm \epsfbox{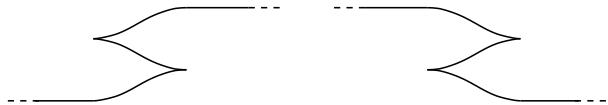}
\endrelabelbox}
  \caption{Legendrian `zigzags'.}
  \label{zigzag} 
\end{figure}

{\bf Contact $r$-surgery with $r<0$.} Let $K$ be the
Legendrian knot along which surgery is to be performed.
Write $r$ as a continued fraction
\[ r_1+1-\cfrac{1}{r_2-\cfrac{1}{\dotsb -\cfrac{1}{r_n} }}\]
with integers $r_1,\ldots ,r_n\leq -2$, cf.~\cite{dige}.
Let $K_1$ be the Legendrian knot represented by the front
projection of $K$ with $|r_1+2|$ additional `zigzags' as in
Figure~\ref{zigzag} (some of which may be of the type
on the left, some of the other type).

For $i=2,\ldots ,n$, let $K_i$ be the Legendrian push-off of
$K_{i-1}$, represented by a parallel copy of the front projection
of $K_{i-1}$ (with the appropriate crossings with the front
projection of~$K_{i-1}$) and with $|r_i+2|$ additional zigzags.

Then a contact $r$-surgery along $K$ corresponds to a sequence
of contact $(-1)$-surgeries along $K_1,\ldots ,K_n$. As observed
in~\cite{dige}, the different choices for the extension of the
contact structure in the process of a contact $r$-surgery
correspond exactly to the different choices of left or right zigzags.

\begin{figure}[h]
\centerline{\relabelbox\small
\epsfxsize 12cm \epsfbox{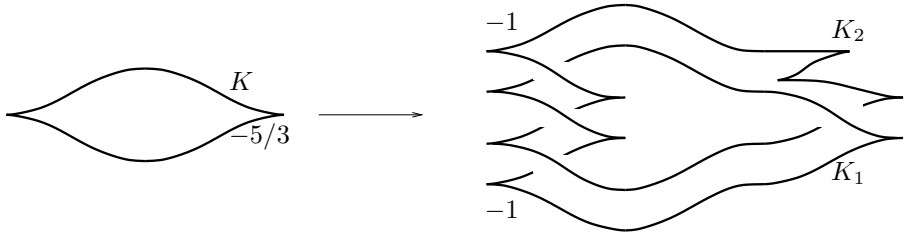}
\extralabel <-1cm,.7cm> {$K_1$}
\extralabel <-5.6cm,.2cm> {$-1$}
\extralabel <-1cm,2.6cm> {$K_2$}
\extralabel <-5.6cm,2.7cm> {$-1$}
\extralabel <-9cm,1.9cm> {$K$}
\extralabel <-9cm,1.2cm> {$-5/3$}
\endrelabelbox}
  \caption{An example for the algorithm.}
  \label{fivethirds} 
\end{figure}

For instance, for $r=-5/3$ we have $r_1=r_2=-3$. Thus, contact
$(-5/3)$-surgery along the Legendrian knot $K$ depicted in
Figure~\ref{fivethirds} is equivalent to a couple
of contact $(-1)$-surgeries along the knots $K_1$, $K_2$.
Here we have to choose an additional
zigzag for $K_1$, and one more for $K_2$. This amounts to four
different possibilities of performing this surgery.

\begin{remark}
In \cite{dige} the sequence of $(-1)$-surgeries replacing a contact
$r$-surgery was defined iteratively, each surgery being performed along
the Legendrian spine of the solid torus glued in when performing the
preceding surgery. There are two ways to see that this is equivalent
to performing successive surgeries along Legendrian push-offs: Assume
$Y'$ is obtained from $(Y,\xi )$ by contact $(-1)$-surgery along
a Legendrian knot~$K$, and write
\[ Y'=(Y\setminus\nu K)\cup S^1\times D^2\]
as before. In the handle picture of \cite[Section~3]{dige}, one can check
that the belt sphere of the $2$-handle corresponding to this surgery
is Legendrian isotopic in $Y'$ to a Legendrian knot $K'\subset Y\setminus
\nu K\subset Y'$ which, when regarded as a knot in~$Y$, is a Legendrian
push-off of~$K$. Alternatively, the Legendrian push-off of a Legendrian knot
$K$ is a knot Legendrian isotopic to $K$ and isotopic on $\partial (\nu K)$
to either of the dividing curves on that convex surface (cf.~\cite{dige01}
for these concepts). The same is true for the spine of the glued in
$S^1\times D^2$, and the gluing is defined by the matching of these
dividing curves.
\end{remark}

{\bf Contact $r$-surgery with $r>0$.} Write $r=p/q$ with $p,q$ coprime
positive integers. Choose a positive integer $k$ such that $q-kp<0$, and
set $r'=p/(q-kp)$. Let $K_1,\ldots ,K_k$ be $k$ successive Legendrian push-offs
of a Legendrian knot~$K$. Then contact $r$-surgery along $K$
is equivalent to contact $(+1)$-surgeries along $K$ and
$K_1,\ldots ,K_{k-1}$, and a contact $r'$-surgery along~$K_k$.
\section{Spin$^c$ structures on 3- and 4-manifolds}
\label{second}
\subsection*{$2$-plane fields and spin$^c$ structures on $3$-manifolds}
In the following we should like to describe surgery diagrams for
contact structures on various 3-manifolds, including all overtwisted 
structures. Since, by~\cite{elia89},
these latter contact structures (up to isotopy)
are in one-to-one correspondence  with oriented 2-plane fields (up to 
homotopy), we begin our discussion by a review of 2-plane fields on
3-manifolds, see~\cite{gomp98} and cf.\ also the discussion in~\cite{geig}
and~\cite{kupe96}.

Let us fix a closed, oriented 
3-manifold $Y$ and consider the space $\Xi (Y)$ of oriented
2-plane fields on $Y$. By considering the oriented normal unit vector
field, we see that the elements of $\Xi (Y)$ are in one-to-one
correspondence with the elements of the space
of vector fields of unit length.

\begin{definition}
{\rm 
Two nowhere vanishing vector fields $v_1$ and $v_2$ are
said to be \emph{homologous} if $v_1$ is homotopic to
$v_2$ outside  a ball $D^3\subset Y$ (through nowhere vanishing
vector fields). An equivalence class of homologous vector fields
is a \emph{spin$^c$ structure}
on $Y$. The set of all spin$^c$ structures is denoted by 
$Spin^c(Y)$.}
\end{definition}

\begin{remark}
{\rm
Traditionally, spin$^c$ structures are defined as lifts of the
orthonormal frame bundle of $Y$ to a principal bundle with
structure group $Spin^c(3)=U(2)$. The equivalence with the definition
given above was observed by Turaev~\cite{tura97}.
}
\end{remark}

Let $\t _{\xi }$ denote the spin$^c$ structure induced
by $\xi \in \Xi (Y)$ (by taking the oriented normal of the
$2$-plane field); this $\t_{\xi}$ depends only on the homotopy
class $[\xi ]$ of~$\xi$. The induced map $[\xi ] \mapsto \t _{\xi }$ will be
denoted by $p\colon\pi_0(\Xi (Y))\to Spin ^c(Y)$; it is obviously
surjective.
It is easy to verify that if $p([\xi _1])=p([\xi _2])$ then we have
equality of first Chern classes $c_1(\xi _1)= c_1(\xi _2)\in H^2(Y)$
(where we regard the oriented $\bfr^2$-bundles $\xi_i$,
uniquely up to homotopy, as complex line bundles).
Therefore we can define the first Chern class of a
spin$^c$ structure $\t \in Spin ^c(Y)$. For the following
standard fact cf.~\cite{tura97}.

\begin{proposition}
The second cohomology group $H^2 (Y; \bfz )$ acts freely and
transitively on $Spin ^c(Y)$. If this action is denoted
by $\t \otimes a$ for $\t \in Spin^c(Y)$ and $a\in H^2(Y;\bfz )$
then $c_1(\t \otimes a)=c_1(\t )+2a$. In particular, if
$H^2(Y; \bfz )$ has no $2$-torsion, then a spin$^c$ structure
$\t$ is uniquely specified by its first Chern class $c_1 (\t )$.
\end{proposition}

For $\t \in
Spin^c(Y)$ the fibre $p^{-1}(\t )$ can be easily identified with the homotopy
classes of $2$-plane fields obtained by taking the connected
sum of $(Y, \xi )$ (where $[\xi ]\in p^{-1}(\t ))$ with the elements of
\[ \{ (S^3, \eta ) \mid \eta {\mbox { is an oriented
$2$-plane field on }} S^3\}\]
(after pasting the $2$-plane fields together). In this way we get a
transitive but not necessarily free ${\bfz}$-action on that fibre.

For $\t \in Spin ^c (Y)$ we denote the divisibility of the
(well-defined) first Chern class $c_1(\t )\in
H^2 (Y; \bfz )$ by $d(\t )$ (which is set to zero if $c_1(\t )$ is torsion).
In the following lemma note that $\bfz_0=\bfz$.

\begin{lemma}[{\cite[Prop.~4.1]{gomp98}}]
The fibre $p^{-1}(\t )\subset \pi _0 (\Xi (Y))$ admits a free and transitive
$\bfz _{d(\t )}$-action. \qed
\end{lemma}

Therefore, for a  spin$^c$ structure whose first Chern class is torsion,
the obstruction to homotopy of two $2$-plane fields both inducing that
given spin$^c$ structure can be captured by a single number.
This obstruction (frequently called the 3-dimensional invariant $d_3$
of $\xi $) can be described as follows:
Suppose that a compact almost complex 4-manifold  $(X,J)$
is given such that
$\partial X=Y$. (Recall that an \emph{almost complex structure} on $X$ 
is a bundle homomorphism $J\colon TX\to TX$ with $J^2=-$id$_{TX}$.) 
The almost complex structure naturally induces a 2-plane field $\xi $ on $Y$
by taking the complex tangencies in $TY$, i.e., $\xi =TY\cap J(TY)$.
Write $\sigma (X),\chi (X)$ for the signature and Euler characteristic
of $X$, respectively.

\begin{theorem}[{\cite[Thm.~4.16]{gomp98}}]\label{d3}
For $c_1(\xi )$ a torsion class, the rational number
\[ d_3(\xi )=\frac{1}{4}\bigl( c_1^2(X,J)-3\sigma (X)-2\chi (X)\bigr)\]
is an invariant of the homotopy type of the
$2$-plane field $\xi$. Moreover, two $2$-plane fields $\xi _1$ and $\xi _2$
with $\t_{\xi_1}=\t_{\xi_2}$ and  $c_1(\t _{\xi _i})=c_1(\xi_i)$
a torsion class are homotopic if and only if $d_3(\xi _1)=d_3(\xi _2).$ \qed
\end{theorem}

\begin{remark}
It is fairly easy to see that for $Y=S^3$ the 3-dimensional invariant
$d_3$ of a 2-plane field lies in $\bfz +\frac{1}{2}$:
for any characteristic vector, hence for $c_1(X,J)$ of an 
almost-complex structure, we have $c_1^2(X, J)\equiv \sigma (X)$ (mod 8)
and $\frac{1}{2}(\sigma (X)+\chi (X))=\frac{1}{2}-b_1(X)+b^+_2(X)$.
The 3-dimensional invariant $d_3$ of $(S^3, \xi _{st})$ (as defined by
Theorem~\ref{d3}) is $-\frac{1}{2}$, since we can regard $(S^3, \xi
_{st})$ as the boundary of the unit disc in $\bfc ^2$.
\end{remark}
\subsection*{Almost complex structures and spin$^c$ structures on
$4$-manifolds}
Let $X$ be a compact 4-manifold, possibly with nonempty boundary $\partial X$.
By a reasoning similar to the $3$-dimensional situation one can see
that an almost complex structure defined on the complement of
finitely many points of $X$ gives rise to a spin$^c$ structure on~$X$.
(This is because both
$S^3$ and $D^4$ admit unique spin$^c$ structures.)
It is fairly easy to see that two such almost complex structures
induce the same spin$^c$ structure if and only if they are homotopic on the
2-skeleton of $X$. This motivates the following definition:

\begin{definition} \label{spinc4}
{\rm 
Two almost complex structures $J_1, J_2$ defined on the complement of
finitely many points in $X$ are \emph{homologous} if there is a compact
$1$-manifold $C\subset X$ containing the finitely many points
where the $J_i$ are undefined such that $J_1$ is homotopic to $J_2$ on $X-C$ 
(through almost complex structures). An equivalence class of homologous almost
complex structures is called a spin$^c$ structure.
The set of spin$^c$ structures on $X$ is denoted by
$Spin^c(X)$.}
\end{definition}

In analogy with the $3$-dimensional case, there is a well-defined notion
of a first Chern class $c_1(\s )$ for $\s\in Spin^c(X)$. The image of the
map $c_1\colon Spin^c(X)\rightarrow H^2(X;{\bfz})$ turns out to equal the
set
\[ \{ c\in H^2(X; \bfz )\mid c\equiv w_2(X)\; \mbox{\rm mod}\; 2\} \]
of {\em characteristic elements}. Once again, $H^2(X;{\bfz})$ acts
freely and transitively on $Spin^c(X)$; we denote this action by
$(\s ,a)\mapsto\s\otimes a$. Again we have $c_1(\s\otimes a)=
c_1(\s )+2a$. Therefore, if $H^2(X;{\bfz})$ has no $2$-torsion, for
instance if $X$ is simply connected, then a spin$^c$ structure $\s$
is uniquely determined by its first Chern class~$c_1(\s )$.

If $Y$ is a $3$-dimensional submanifold of~$X$, then a spin$^c$ structure
on $X$ naturally induces a spin$^c$ structure on $Y$ by taking the orthogonals
of the complex tangencies in $TY$.
\subsection*{Homological data of $2$-handlebodies.}
In our later arguments we shall make computations involving homology
and cohomology classes on 2-handle\-bodies and on their boundaries. So
let us assume that the 4-manifold $X$ is given by the framed link $\Li
=((K_1, n_1), \ldots , (K_t, n_t))\subset S^3$, i.e., we attach copies
of $D^2\times D^2$ along $\partial D^2 \times D^2$ to $D^4$ along $\nu
K_i\subset \partial D^4=S^3$ with the specified framing $n_i$.  (For
more about such \emph{Kirby diagrams} see \cite{gost99}.  Note that we
only deal with the case when $X$ is decomposed into one 0-handle and a certain
number $t$ of 2-handles.)

Obviously $\pi _1 (X)=1$, and $H_2 (X; \bfz )$ is
generated by the fundamental classes $[\Sigma_i]$ of the surfaces
$\Sigma _i$ we get by gluing a Seifert surface $F_i$ of $K_i$ to the
core disc of the $i^{th}$ handle. The
intersection form in this basis of $H_2(X; \bfz )$ is simply the
linking matrix of $\Li$, with the framing coefficients $n_i$ in the
diagonal.

Let $N_i$ denote a small normal disc to $K_i$ in $S^3$ and
$\mu _i =\partial N_i$.  An orientation on the knot $K_i$ will give an
orientation of $\Sigma _i$ (by requiring that the orientation
of $K_i$ be the boundary orientation of the Seifert surface~$F_i$).
Together with the orientation of the ambient 3-manifold $S^3$, the
orientation of $K_i$ will induce an orientation on $N_i$ as
well. We can then give $\mu_i=\partial N_i$ the boundary orientation. In
the knot diagrams below the
orientation of $K_i$ will be denoted by a little arrow
next to the diagram of the knot.

It is easy to see that the relative
homology classes $[N_i]$ freely generate $H_2(X, \partial X; \bfz )$,
while $H_1(\partial X; \bfz )$ is generated by the homology
classes $[\mu_i]$ of the circles $\mu_i=\partial N_i$ ($i=1, \ldots , t$).
The long exact sequence of the
pair $(X, \partial X)$ reduces to
$$0\to H_2(\partial X; \bfz )\to H_2(X; \bfz )
\stackrel{\varphi _1}{\longrightarrow}
H_2(X, \partial X; \bfz )
\stackrel{\varphi _2}{\longrightarrow}
H_1(\partial X; \bfz )\to 0,$$
since the condition $\pi _1(X)=1$ implies
\[ H_1(X; \bfz )=0=H^1(X; \bfz )\cong H_3(X,\partial X; \bfz ). \]
The maps $\varphi _1$ and $\varphi _2$ are easy to describe 
in the above bases:
With $\ell k(K_i, K_j)$ denoting the linking number of $K_i$ and
$K_j$ for $i\neq j$ and $\ell k(K_i, K_i)= n_i$ we have
$$\varphi _1 ([\Sigma _i])=\sum _{j=1}^t  \ell k (K_i, K_j)[N_j] ;$$
furthermore
$$\varphi _2 ([N_i])=[\mu _i].$$
(For details of the argument see \cite{gost99}.)
For a cohomology class $c\in H^2 (X; \bfz )$ denote by
$c(\Sigma_i)=\langle c, [\Sigma _i]\rangle\in\bfz$ its evaluation
on~$\Sigma_i$. Then the Poincar\'e dual
$PD (c)\in H_2(X, \partial X; \bfz )$ is equal to $\sum _{i=1}^t
c(\Sigma _i)[N_i]$. The image $\varphi _2(PD (c))$ gives
a description of $PD (c\vert _{\partial X})$ in terms of the
1-homologies $[\mu _i]$. Exactness of the sequence implies that the
relations among the $[\mu _i]$ are simply given by the
expressions $\varphi _1([\Sigma _i])$ with $[N_i]$ substituted by 
$[\mu _i]$. These relations help to simplify $PD (c\vert _{\partial X})$.
If that class is a torsion element then for appropriate $n\in \bfn $ the
class $PD(n\cdot c)$ maps to zero under $\varphi _2$, hence it is 
the image of a class $C\in H_2(X;\bfz )$
under $\varphi _1$. In that case we can compute 
$c^2$ as $c^2=C^2/n^2$.
\section{Computation of homotopy invariants of contact structures}
\label{third}
From a surgery presentation of $(Y,\xi )$ we now wish to determine some
homotopy invariants of $\xi $. 
The surgery diagram can be considered as a Kirby diagram for a 
4-manifold $X$ with boundary $Y$. Consider $(S^3, \xi _{st})$ as the
boundary of the standard disc $D^4\subset \bfc ^2$, equipped with its
standard (almost) complex structure.

\begin{proposition}[\cite{elia90}, {\cite[Prop.~2.3]{gomp98}}]
\label{prop:attach-1}
If a $2$-handle $H$ is attached along a Legendrian knot $K\subset 
(S^3, \xi _{st})$ with framing $(-1)$  (i.e.\ one left
twist added to the contact framing) then the above standard
complex structure extends as an (almost) complex structure $J$ to
$D^4\cup H$ inducing the surgered contact structure on the boundary.
Moreover, $c_1(D^4\cup H, J)$ evaluates on the homology class given by
$K$ (in the sense of the previous section)
as ${\rm {rot}}(K)$.\qed
\end{proposition} 

\begin{remark}
In fact, Eliashberg~\cite{elia90}
showed that the Stein structure of $D^4$ extends as
a Stein structure to $D^4\cup H$, cf.~\cite{foko}.
\end{remark}

We now want to study the related question for contact $(+1)$-surgeries. Thus,
let $X=D^4\cup H$ be the handlebody corresponding to a contact
$(+1)$-surgery on a Legendrian knot $K\subset (S^3,\xi_{st})=\partial
D^4$. The contact structure $\xi$ on $\partial X$ determined by the
surgery defines an almost complex structure $J$ (on~$X$) along~$\partial X$,
unique up to homotopy: require, firstly, $\xi$ to be $J$-invariant (and
the orientation of $\xi$ induced by $J$ to coincide with the given one)
and, secondly, $J$ to map the outward normal along $\partial X$ to a vector
positively transverse to~$\xi$.

That $J$ extends to the complement of a $4$-disc
$D_H\subset\mbox{\rm int}(H)\subset X$, for there is no obstruction to
extending $J$ over the cocore $2$-disc of the $2$-handle, and
$X-D_H$ deformation retracts onto the union of $\partial X$
and that cocore disc. In particular, there is a class $c\in H^2(X;{\bfz})$
that restricts to $c_1(\xi )=c_1(J)$ on $\partial X$, and whose
mod~$2$ reduction equals~$w_2(X)$; the existence of such a class
(which conversely implies the existence of $J$ on $X-D_H$) can also be shown
by a purely homological argument.

Let $\xi_H$ be the plane field on $\partial D_H=S^3$ induced by~$J$,
where $\partial D_H$ is given the orientation as boundary
of $D_H\subset X$ rather than the boundary orientation of $\partial
(X-D_H)$. By \cite{gomp98},
cf.~\cite[Thm.~11.3.4]{gost99} and the discussion preceding it,
there is an almost complex manifold $(W,J_W)$ with $\partial W=S^3$ such
that $J_W$ induces the plane field $\xi_H$ on the boundary. With the help
of $W$ one can compute the invariant $d_3(\xi_H)$. Moreover, by the
proof of that same quoted theorem the $d_3$-invariant behaves
additively in the sense that
\begin{eqnarray*}
d_3(\xi ) & = & \frac{1}{4}\bigl( c_1^2(X-D_H,J)-3\sigma (X-D_H)
                -2\chi (X-D_H)\bigr) +d_3(\xi_H)\\
          & = & \frac{1}{4}\bigl( c^2-3\sigma (X)-2\chi (X)\bigr)
                  +d_3(\xi_H)+\frac{1}{2}.
\end{eqnarray*}

\begin{remark}
\label{rem:pm}
{\rm
As part of the following proposition we shall see that $d_3(\xi_H)=1/2$.
This equals the $d_3$-invariant of the standard contact structure on~$S^3$,
regarded as the boundary of ${\bfc}P^2-D^4$ (i.e.\ with the
opposite of the usual
orientation, which causes the sign change of the $d_3$-invariant,
cf.~\cite[Thm.~11.3.4]{gost99} again). Thus an equivalent way of phrasing
the result $d_3(\xi_H)=1/2$ is that the almost complex structure
defined near $\partial X$ extends over $X\# {\bfc}P^2$, coinciding with the
standard structure near the $2$-skeleton of~${\bfc}P^2$. (See also
Section~\ref{section:+1} below.)

On the other hand,
recall from \cite{dige} that contact $(+1)$-surgery can be
regarded as a symplectic handlebody surgery on the {\em concave} end
of a symplectic cobordism. In particular, we may regard
$(S^3,\xi_{st})$ ({\em with reversed orientation})
as the concave boundary of ${\mathbb C}P^2-D^4$ with
its standard K\"ahler structure, and contact $(+1)$-surgery along $K$
corresponds to adding a symplectic $2$-handle to ${\mathbb C}P^2-D^4$
along its boundary. This implies that the contact structure
on $\partial X$ {\em with reversed orientation} is induced from an
almost complex structure on $\overline{X}\# {\bfc}P^2$, again
coinciding with the standard structure near the $2$-skeleton of~${\bfc}P^2$.
(Here $\overline{X}$ denotes $X$ with reversed orientation.)

Thus, we can glue $X\#{\bfc}P^2$ and $\overline{X}\#{\bfc}P^2$ along
their common boundary (with opposite orientations) to obtain an almost
complex manifold
\[ {\bfc}P^2\# X\cup\overline{X}\#{\bfc}P^2={\bfc}P^2\# DX\#{\bfc}P^2,\]
where $DX$ denotes the double of~$X$, which in the present situation
is diffeomorphic to $S^2\times S^2$ or ${\bfc}P^2\#\overline{{\bfc}P^2}$,
cf.~\cite[Cor.~5.1.6]{gost99}. Indeed, a homological calculation similar
to the following proof shows that ${\bfc}P^2\# DX\#{\bfc}P^2$ admits an
almost complex structure, standard near the $2$-skeleta of the
${\bfc}P^2$-summands, which splits in the way described.

If $X$ is a handlebody corresponding to $n$ contact $(+1)$-surgeries,
then the contact manifold $\partial X$ is boundary of the almost complex
manifold $X\# n{\bfc}P^2$; with reversed orientation it is the
boundary of $\overline{X}\#{\bfc}P^2$. Again one checks that
$n{\bfc}P^2\# DX\#{\bfc}P^2$ admits an appropriate almost complex structure.
($DX$ is diffeomorphic to $nS^2\times S^2$ or $n{\bfc}P^2\# n
\overline{{\bfc}P^2}$, cf.~~\cite[Cor.~5.1.6]{gost99}.)
}
\end{remark}

\begin{proposition}
\label{prop:+1}
Let  $K\subset (S^3,\xi_{st})$ be a Legendrian knot with $\tb (K)\neq 0$. 
If the handlebody $X$ is obtained by attaching
a $2$-handle $H$ to $D^4$ along $K$ with
framing $(+1)$ (one right twist added to the contact framing), then the almost
complex structure defined near $\partial X$ extends over
$X-D_H$, in the previously introduced notation, such that
$d_3(\xi_H)=1/2$.  Moreover, the corresponding
class $c\in H^2(X;{\bfz})$ evaluates on the homology class given
by $K$ as $\rot (K)$.
\end{proposition}

\begin{proof}
Consider Legendrian push-offs $K_1,\ldots ,K_n,K_1',\ldots, K_n'$
of $K$ (it would be enough to study the cases $n=1$ or~$2$). Perform
contact $(+1)$-surgeries on $K_1,\ldots ,K_n$ and contact $(-1)$-surgeries
on $K_1',\ldots ,K_n'$. By Lemma~\ref{lem:cancel} the
resulting manifold is $(S^3,\xi_{st})$.

Let
\[ \Sigma_1,\ldots ,\Sigma_n,\Sigma_1',\ldots ,\Sigma_n'\]
be the corresponding surfaces in
\[ X=X_n=D^4\cup H_1\cup\ldots\cup H_n\cup H_1'\cup\ldots\cup H_n'\]
in the notation of the preceding section. Write
$c=c_{(n)}\in H^2(X_n;{\mathbb Z})$
for the class defined by the almost complex structure on $X$ with
discs $D_{H_1},\ldots ,D_{H_n}$
removed. By Proposition~\ref{prop:attach-1} we have
$c(\Sigma_i')=\rot (K)$, $i=1,\ldots ,n$. Set $k=c(\Sigma_i)$. Then, again
by the preceding section (and in the notation used there),
\[ PD(c)=k\sum_{i=1}^n [N_i]+\rot (K)\sum_{i=1}^n [N_i'].\]
This can be written as $PD(c)=\varphi_1(C)$ with a unique class
$C\in H_2(X;{\mathbb Z})$ (since $H_1(\partial X)=H_2(\partial X)=0$).
We have
\[ \varphi_1([\Sigma_i])=\tb (K)\sum_{j=1}^n\bigl( [N_j]+[N_j']\bigr)
+[N_i]\]
and
\[ \varphi_1([\Sigma_i'])=\tb (K)\sum_{j=1}^n\bigl( [N_j]+[N_j']\bigr)
-[N_i'].\]
Write
\[ C=\sum_{i=1}^n\bigl( a_i[\Sigma_i]+a_i'[\Sigma_i']\bigr) .\]
Then the coefficients $a_i,a_i'$ are found as solutions of
the linear equation
\[ M_{\mbox{\scriptsize\rm tb} (K)}\left( \begin{array}{c}a_1\\
\vdots\\ a_n\\a_1'\\
\vdots\\ a_n'\end{array}\right)=
\left( \begin{array}{c}k\\ \vdots\\ k\\ \rot (K)\\ \vdots\\ \rot(K)
\end{array}\right) ,\]
where $M_{\mbox{\scriptsize\rm tb} (K)}$ is the matrix
\[ M_{\mbox{\scriptsize\rm tb} (K)} =
\tb (K) E_{2n}+ \left( \begin{array}{cc}I_n&0\\0&-I_n
\end{array}\right),\]
with $E_{2n}$ the $(2n\times 2n)$-matrix having all entries equal to~$1$,
and $I_n$ the $(n\times n)$ unit matrix.

It follows that
\[ a_1=\cdots = a_n=k-n\bigl( k-\rot (K)\bigr) \tb (K)\]
and
\[ a_1'=\cdots = a_n'=-\rot (K)+n\bigl( k-\rot (K)\bigr) \tb (K),\]
whence
\begin{eqnarray*}
c^2\; =\; C^2 & = & \bigl( a_1,\ldots ,a_n ,a_1',\ldots ,a_n'\bigr)
M_{\mbox{\scriptsize\rm tb} (K)}\left( \begin{array}{c}a_1\\ \vdots\\ a_n\\
a_1'\\ \vdots \\ a_n'\end{array}\right) \\
  & = & \bigl( a_1,\ldots ,a_n ,a_1',\ldots ,a_n'\bigr)
\left( \begin{array}{c}k\\ \vdots\\ k\\ \rot (K)\\ \vdots\\ \rot(K)
\end{array}\right)\\
  & = & n\bigl( k^2-\rot^2(K)\bigr)-n^2\tb (K)\bigl( k-\rot (K)\bigr)^2.
\end{eqnarray*}

The signature of $M_{\mbox{\scriptsize\rm tb} (K)}$ (which is the same
as the signature of~$X$)
is equal to zero (this follows from the fact
that $M_{\mbox{\scriptsize\rm tb} (K)}$
remains nonsingular if $\tb (K)$ is replaced by
any real parameter). The Euler characteristic of $X$ is $1+2n$.

From the discussion preceding the present proposition we deduce
\begin{eqnarray*}
\lefteqn{-\frac{1}{2}=d_3(S^3,\xi_{st})  = }\\
 & & = \frac{1}{4}\bigl( c_1^2(X-\bigcup_i D_{H_i},J)-3\sigma (X)
-2\chi (X)\bigr) +n\bigl( d_3(\xi_H)+\frac{1}{2}\bigr)\\
 & & = \frac{1}{4}\bigl[ n\bigr( k^2-\rot^2(K)\bigr) -n^2\tb (K)
\bigl( k-\rot (K)\bigr)^2\bigr] +\\
 & & \;\;\;\;\;\;\;\mbox{}+n\bigl( d_3(\xi_H)-\frac{1}{2}\bigr)
      -\frac{1}{2}.
\end{eqnarray*}
This is true for any $n\in{\mathbb N}$, from which we conclude, for
$\tb (K)\neq 0$, that $k=\rot (K)$ and $d_3(\xi_H)=1/2$.
\end{proof}

\begin{remark}
\label{rem:tb}
{\rm
The result $d_3(\xi_H)=1/2$ remains true even if $\tb (K)=0$. This can
be seen from the description of contact $(+1)$-surgery in~\cite{dige}
as a symplectic handlebody surgery on the concave end of a symplectic
cobordism. Indeed, this description provides a unique model for
contact $(+1)$-surgery, so that the obstruction for extending
the almost complex structure over the handle is independent of~$\tb (K)$.

In the case $\tb (K)=0$
the above argument only yields $k=\pm {\rm {rot}}(K)$.
The quickest way to see that $k={\rm {rot}}(K)$ in this case as well
is the following: Since, as just remarked, contact $(+1)$-surgery also admits
a handlebody description, one can mimic the argument
of~\cite[Prop.~2.3]{gomp98}, where the corresponding result was shown for
contact $(-1)$-surgeries. Checking all the relevant signs might be
tedious, but again the argument shows that $k$ does not depend on $\tb (K)$,
so our result $k=\rot (K)$ for $\tb (K)\neq 0$ in fact also holds
in the case $\tb (K)=0$.

Since we shall not use the result for ${\rm {tb}}(K)=0$ in our subsequent
arguments, we defer an explicit discussion of these issues to
Section~\ref{section:+1}.
}
\end{remark}

\begin{corollary}
\label{cor:d3}
Suppose that $(Y, \xi )=\partial X$, with $c_1(\xi )$ torsion,
is given by contact $(\pm 1)$-surgery
on a Legendrian link ${\mathbb {L}}\subset (S^3, \xi _{st})$
with $\tb (K)\neq 0$ for each $K\subset {\mathbb L}$ on which we
perform contact $(+1)$-surgery. Then
$$d_3(\xi )=\frac{1}{4}\bigl( c^2-3\sigma (X)-2\chi (X)\bigr) +q,$$
where $q$ denotes the number of components in ${\mathbb {L}}$ on which we
perform $(+1)$-surgery, and $c\in H^2 (X; \bfz )$ is the cohomology class
determined by $c(\Sigma_K)={\rm {rot}}(K)$ for each $K\subset \Li$. Here
$[\Sigma_K]$ is the homology class in $H_2(X)$ determined by $K\subset
S^3$ (i.e.\ Seifert surface of $K$ glued with core disc of corresponding
handle).
\end{corollary}

\begin{proof}
The contact manifold $(Y,\xi )$ is the boundary of the almost complex
manifold $X\# q{\mathbb C}P^2$ (such that $\xi$ is given by the complex
tangencies in $Y=\partial X$), with first Chern class
\[ c_1=c+(3,\ldots ,3)\in H^2(X;{\mathbb Z})\oplus
qH^2({\mathbb C}P^2;{\mathbb Z}),\]
which satisfies $c_1^2=c^2+9q$.

Moreover,
\[ \sigma (X\# q{\mathbb C}P^2)=\sigma (X)+q \]
and
\[ \chi (X\# q{\mathbb C}P^2)=\chi (X)+q.\]
Hence
\begin{eqnarray*}
d_3(\xi) & = & \frac{1}{4}\bigl( c^2+9q-3(\sigma (X)+q)-2(\chi (X)+q) \bigr) \\
         & = & \frac{1}{4}\bigl( c^2-3\sigma (X)-2\chi (X)\bigr) +q.
\end{eqnarray*}
\end{proof}

\section{Surgery diagrams for overtwisted contact $3$-manifolds}
\label{fourth}
Our next goal is to draw surgery diagrams for all overtwisted 
contact structures on a given 3-manifold $Y$.
Recall from \cite{elia89} that overtwisted contact structures (up to isotopy)
are in one-to-one correspondence with elements of 
$\pi _0 (\Xi (Y))$. Therefore, in order to find
all the necessary diagrams, we need to find, for each spin$^c$ structure
on~$Y$, a surgery diagram for a contact structure inducing that
spin$^c$ structure, and diagrams for
all overtwisted contact structures on $S^3$. By taking connected sums of
these structures -- which is reflected simply as disjoint union in
the diagrams -- we get all the pictures we wanted. First we 
show how to draw surgery diagrams for all contact structures 
on $S^3$. Then we do the same for $S^1\times S^2$, and finally  we turn 
to the general case.

Notice also that if $(Y, \xi )$ is a contact structure with
$c_1(\t _{\xi})$ torsion and $(S^3, \xi _i )$ is an arbitrary
contact structure, then 
\[ d_3(Y, \xi \# \xi _i )=d_3(Y, \xi )+d_3(S^3 , \xi _i)+\frac{1}{2} .\]
This follows from the fact that under the boundary connected sum
$X\natural X'$ of $4$-manifolds $X,X'$, the signature $\sigma$ and
the number $c^2$ behave additively, whereas
\[ \chi (X\natural X')=\chi (X)+\chi (X')-1.\]

\subsection*{Contact structures on $S^3$}
By Eliashberg's classification~\cite{elia89}, \cite{elia92} we know that
$S^3$ admits a unique tight contact 
structure $\xi _{st}$ (which can be represented by the empty diagram in 
$(S^3, \xi _{st})$), and a unique overtwisted one (up to isotopy)
in each homotopy class of 2-plane fields. Obviously, all these
structures have zero first Chern class; the overtwisted ones can be
distinguished by their 3-dimensional invariant~$d_3$.

\begin{figure}[h]
\centerline{\relabelbox\small
\epsfxsize 12cm \epsfbox{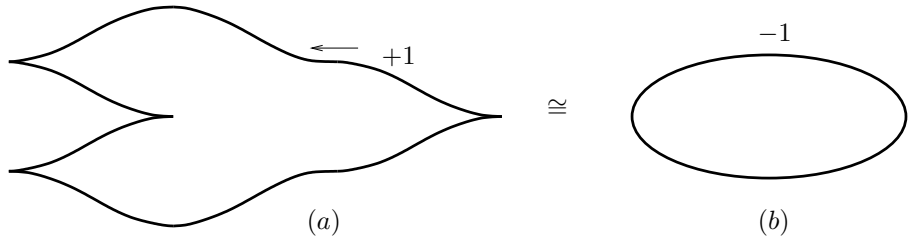}
\extralabel <-7cm,2.2cm> {$+1$}
\extralabel <-2cm,2.5cm> {$-1$}
\extralabel <-4.8cm,1.5cm> {$\cong$}
\extralabel <-8cm,0cm> {$(a)$}
\extralabel <-2cm,0cm> {$(b)$}
\endrelabelbox}
  \caption{Contact structure $\xi _1$ on $S^3$ with $d_3(\xi _1)=1/2$.}
  \label{shark2} 
\end{figure}

\begin{figure}[h]
\centerline{\relabelbox\small
\epsfxsize 12cm \epsfbox{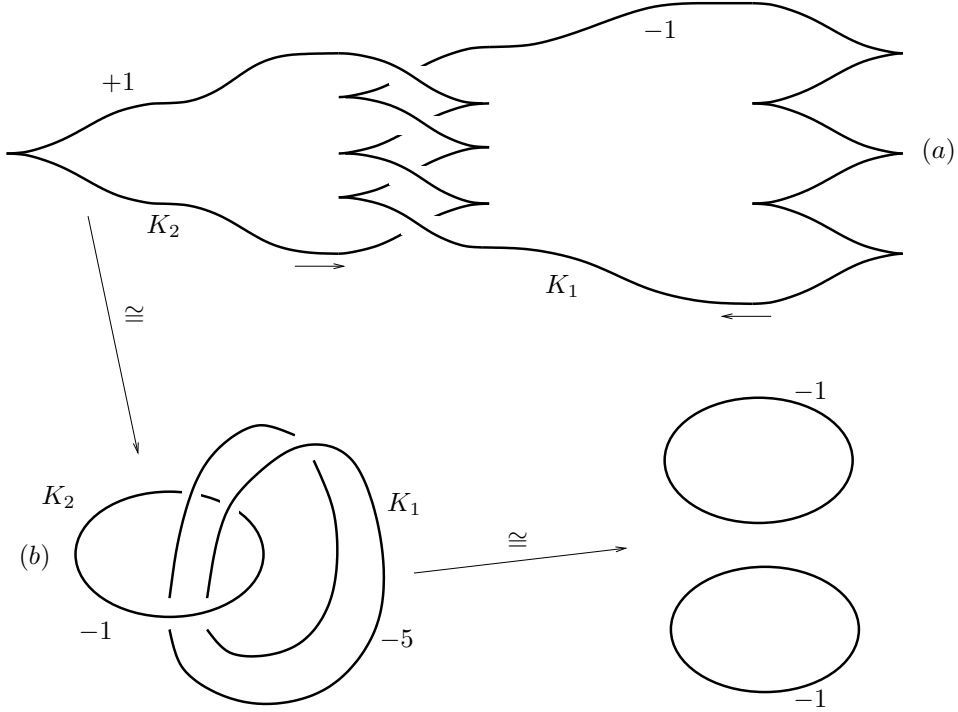}
\extralabel <-1.5cm,.cm> {$-1$}
\extralabel <-1.5cm,4.1cm> {$-1$}
\extralabel <-5.3cm,2.1cm> {$\cong$}
\extralabel <-6.9cm,2.6cm> {$K_1$}
\extralabel <-7.0cm,.8cm> {$-5$}
\extralabel <-11cm,.9cm> {$-1$}
\extralabel <-11.5cm,2.7cm> {$K_2$}
\extralabel <-10.7cm,8.2cm> {$+1$}
\extralabel <-4.8cm,5.5cm> {$K_1$}
\extralabel <-3.5cm,8.9cm> {$-1$}
\extralabel <-10.1cm,6.3cm> {$K_2$}
\extralabel <-10.4cm,5.1cm> {$\cong$}
\extralabel <.2cm,7.3cm> {$(a)$}
\extralabel <-11.8cm,1.9cm> {$(b)$}
\endrelabelbox}
  \caption{Contact structure $\xi_{-1}$ on $S^3$ with
$d_3(\xi_{-1})=-3/2$.}
  \label{negative} 
\end{figure}

\begin{lemma}
\label{lemma:shark}
The surgery diagram of Figure~\ref{shark2}(a)
 gives a contact structure $\xi _1$ 
on  $S^3$ with $d_3(\xi _1)=\frac{1}{2}$.
The surgery diagram of Figure~\ref{negative}(a) gives a contact structure 
$\xi _{-1}$ on  $S^3$ with $d_3(\xi _{-1})=-\frac{3}{2}$.
\end{lemma}

\begin{proof}
By turning the diagrams into smooth surgery diagrams 
(i.e., disregarding the Legendrian position of the surgery
curves and thus the induced contact structure on the result)
and reading the framings not relative to the contact framing, but
relative to the framings
induced by the Seifert surfaces in 
$S^3$, we see that topologically the two surgeries yield~$S^3$.
The equivalence between the surgery descriptions in
Figure~\ref{negative}(b) (even as Kirby diagrams of a $4$-manifold)
is given by a handle slide; cf.~\cite[p.~150]{gost99}.

Here is the computation of the $d_3$-invariants (with notation as
above):

Recall from \cite{gomp98}, \cite{gost99} that for a Legendrian knot~$K$,
represented by its front projection, we have
\[ \tb (K)=\mbox{\rm writhe}(K)-\frac{1}{2}\# (\mbox{\rm cusps})\]
and
\[ \rot (K)=\frac{1}{2}(\# (\mbox{\rm down-cusps})-\#
(\mbox{\rm up-cusps})).\]
Thus, in the first case we have, with the indicated orientation of the
Legendrian knot~$K$, that $\rot (K)=1$.
Hence
\[ PD(c)=c(\Sigma )[N]=\rot (K)[N]=[N].\]
Since the topological framing of $K$ (i.e.\ the framing relative to the
surface framing) is $k=-1$, we have $\varphi_1([\Sigma ])=-[N]$. Therefore
$C=-[\Sigma ]$ and $c^2=C^2=k=-1$. Moreover, the corresponding handlebody
$X=D^4\cup H$ has $\sigma (X)=\mbox{sign}(k)=-1$ and $\chi (X)=2$. Thus,
by Corollary~\ref{cor:d3},
\[ d_3(\xi_1)=\frac{1}{4}(-1+3-4)+1=\frac{1}{2}.\]

In the second case, again with the indicated orientations, we have
\[ \tb (K_1)=-4,\;\; \rot (K_1)=1,\;\; \tb (K_2)=-2,\;\; \rot (K_2)=-1.\]
Furthermore, the linking number $\ell k(K_1,K_2)$ equals $-2$, so the linking
matrix, which describes the homomorphism~$\varphi_1$,
is $\left( \begin{array}{rr}-5&-2\\-2&-1\end{array}\right)$.
With
\[ PD (c)=\rot (K_1)[N_1]+\rot (K_2)[N_2]=[N_1]-[N_2]\]
we find that
the solution of $\varphi_1(C)=PD(c)$ is $C=-3[\Sigma_1]+7[\Sigma_2]$. Thus
\[ c^2=C^2=(-3,7)\left( \begin{array}{rr}-5&-2\\-2&-1\end{array}\right)
\left(\begin{array}{c}-3\\7\end{array}\right)=-10.\]
Moreover, the corresponding handlebody $X=D^4\cup H_1\cup H_2$ has
$\chi (X)=3$ and $\sigma (X)=-2$ (which is obvious from the smooth
surgery description). We conclude
\[ d_3(\xi_{-1})=\frac{1}{4}(-10+6-6)+1=-\frac{3}{2}.\qed\]
\renewcommand{\qed}{}\end{proof}

Using the connected sum operation on the two basic contact structures
$\xi_1$ and $\xi_{-1}$,
we can now draw diagrams for all overtwisted contact structures 
$\xi_i$ on $S^3$ with $d_3(\xi _i)=i-1/2$ ($i \in \bfz )$.
Of course, this procedure will not necessarily provide the most
``economic '' surgery diagram of $\xi _i$.

\vspace{2mm}

Here is a brief sketch of an alternative construction: Let $K_1$ be
a Legendrian knot in $(S^3,\xi_{st})$. Let $K_2$ be the Legendrian
knot obtained from a Legendrian push-off of $K_1$ by adding two
zigzags to its front projection, and perform contact $(+1)$-surgery
on both knots. Topologically, contact $(+1)$-surgery on $K_2$ is
the same as contact $(-1)$-surgery along a Legendrian push-off of
$K_1$, so the resulting manifold is again $S^3$ by
Lemma~\ref{lem:cancel}. Write $\xi$ for the contact structure
on $S^3$ obtained via that surgery.

Equip $K_1$ with an orientation.
By a computation as in the proof of the preceding lemma, one finds
that
if $K_2$ is obtained from a Legendrian push-off of $K_1$
by adding two down-zigzags to its front projection,
then $d_3(\xi )= -\tb (K_1)-\rot (K_1)-1/2$.

Any odd (but no even)
integer can be realised as $\tb (K_1)+\rot (K_1)$ for
a suitable Legendrian knot $K_1$. We leave it as an exercise to the
reader to construct such $K_1$ (see the examples in~\cite{gomp98}
and~\cite{etho01}); that even integers are excluded follows from
\cite[Prop.~2.3.1]{elia93}. Therefore, any overtwisted contact structure
on $S^3$ can be obtained by contact $(+1)$-surgeries on either two
or three Legendrian knots (to realise $d_3=2m-1/2$, $m\in\bfz$, construct a
contact structure $\xi$ on $S^3$ with $d_3(\xi )=(2m-1)-1/2$ by two
$(+1)$-surgeries as just described, then take the connected sum
with $(S^3,\xi_1)$).
\subsection*{Contact structures on $S^1\times S^2$}
According to a folklore theorem of Eliashberg, $S^1\times S^2$ admits a
unique tight contact structure (for a sketch proof see Exercise~6.10
in~\cite{etny01}).

\begin{figure}[h]
\centerline{\relabelbox\small
\epsfxsize 8cm \epsfbox{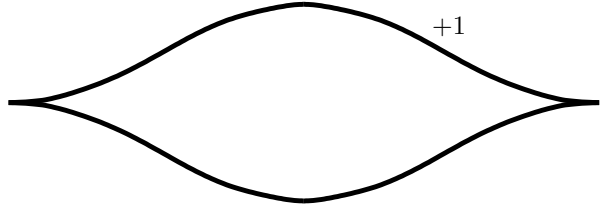}
\extralabel <-2.3cm,2.3cm> {$+1$}
\endrelabelbox}
  \caption{Legendrian unknot producing tight $S^1\times S^2$.}
  \label{saucer} 
\end{figure}

\begin{lemma}
Contact $(+1)$-surgery on the Legendrian unknot (see Figure~\ref{saucer}) 
yields the tight contact structure on $S^1\times S^2$.
\end{lemma}

\begin{proof}
The Legendrian unknot shown in Figure~\ref{saucer} has Thurston-Bennequin
invariant~$-1$, thus contact $(+1)$-surgery corresponds to a
topological $0$-surgery, which produces the manifold $S^1\times S^2$.

For the contact-geometric part of the
proof we use the language of convex surfaces and dividing curves;
for a brief introduction see~\cite{etny01}. By~\cite[Thm.~8.2]{kand97}
and~\cite[Prop.~4.3]{hond00}, for any $k\in{\mathbb Z}$
there is a unique tight contact structure
on $S^1\times D^2$ with a fixed convex boundary with dividing set
consisting of two curves of slope $1/k$, where the meridian corresponds to
slope zero and the longitude $S^1\times \{ p\}$, $p\in\partial D^2$,
to slope~$\infty$. Notice that different values of $k$ simply correspond
to a different choice of longitude. It therefore suffices to show that both
the standard tight contact structure on $S^1\times S^2$ and the
contact structure obtained by the described surgery can be split along
an embedded convex $T^2$ with dividing set as described.

The standard tight contact structure on $S^1\times S^2\subset S^1\times
{\mathbb R}^3$ is given, in obvious notation, by 
\[ \alpha :=x\, d\theta+y\, dz -z\, dy =0.\]
Embed $T^2$ as follows:
\[
\begin{array}{rcl}
T^2 & \longrightarrow & S^1\times S^2\\
(\theta ,\varphi ) & \longmapsto & (\theta ,f(\varphi ),
                     \sqrt{1-f^2(\varphi )}\cos\varphi ,
                     \sqrt{1-f^2(\varphi )}\sin\varphi)
\end{array}\]
with $f (\varphi )=\varepsilon\sin\varphi$ for some $\varepsilon\in (0,1)$.
The tangent spaces of this embedded $T^2$ are spanned by $\partial_{\theta}$
and
\[ v=(0,f',-\frac{ff'}{\sqrt{1-f^2}}\cos\varphi-\sqrt{1-f^2}\sin\varphi ,
           -\frac{ff'}{\sqrt{1-f^2}}\sin\varphi+\sqrt{1-f^2}\cos\varphi ).\]
From $\alpha (\partial_{\theta})=f$ and $\alpha (v)=1-f^2$ we conclude
that the characteristic foliation on $T^2$ is given by
${\displaystyle\partial_{\theta}-\frac{f}{1-f^2}\partial_{\varphi}}$,
which admits the dividing curves $\{\varphi =\pi /2\}$ and $\{\varphi =3\pi
/2\}$.  This means that $T^2$ is a convex torus with dividing set consisting
of two longitudes, as desired.

Now to the same question for the contact structure on $S^1\times S^2$ obtained
via the indicated surgery. First of all, we recall that in the unique
local contact geometric model for the tubular neighbourhood of a Legendrian
knot, the boundary of that neighbourhood is a convex torus with
dividing set consisting of two copies of the longitude
determined by the contact framing, cf.~\cite{dige}. Write $K$ for the
Legendrian knot of Figure~\ref{saucer} and $\nu K$ for a (closed) tubular
neighbourhood. Further, we denote the meridian of $\partial (\nu K)$ by
$\mu$, and by $\lambda$ the longitude determined by $\ell k (\lambda ,K)=0$.

Then $S^3-\mbox{\rm int}(\nu K)$ is a solid torus with meridian
$\overline{\mu}=\lambda$ and a longitude $\overline{\lambda}=\mu$.
Since $\tb (K)=-1$, the longitude $\lambda_c$ determined by the
contact framing is
\[ \lambda_c=\lambda -\mu =\overline{\mu}-\overline{\lambda},\]
which is a longitude of $S^3-\mbox{\rm int}(\nu K)$, so the tight contact
structure on that piece has a convex boundary of the kind described above.

The surgered manifold ($(+1)$-surgery with respect to the framing
given by~$\lambda_c$) is given by
\[ (S^3-\nu K)\cup N_0,\]
where $N_0$ is a solid torus, with meridian $\mu_0$ and longitude
$\lambda_0$ of $\partial N_0$ being glued to
$\partial (\nu K)$ by
\[ \mu_0\longmapsto -\lambda_c-\mu =-\lambda ,\;\;\;
   \lambda_0\longmapsto\mu .\]
Observe that the curve $-\mu_0-\lambda_0$ is glued to a dividing
curve $\lambda_c=\lambda -\mu$. So the extension of the contact
structure over $N_0$ in the process of contact surgery is given by
the unique tight contact structure with convex boundary having
two copies of the longitude $-\mu_0-\lambda_0$ as dividing set.
This concludes the proof.
\end{proof}

\begin{remark}
{\rm
An alternative proof of this lemma, deducing tightness from
the non-vanishing of the corresponding Heegaard-Floer invariant,
is given in~\cite[Lemma~4]{listb}.
}
\end{remark}

\begin{figure}[h]
\centerline{\relabelbox\small
\epsfxsize 12cm \epsfbox{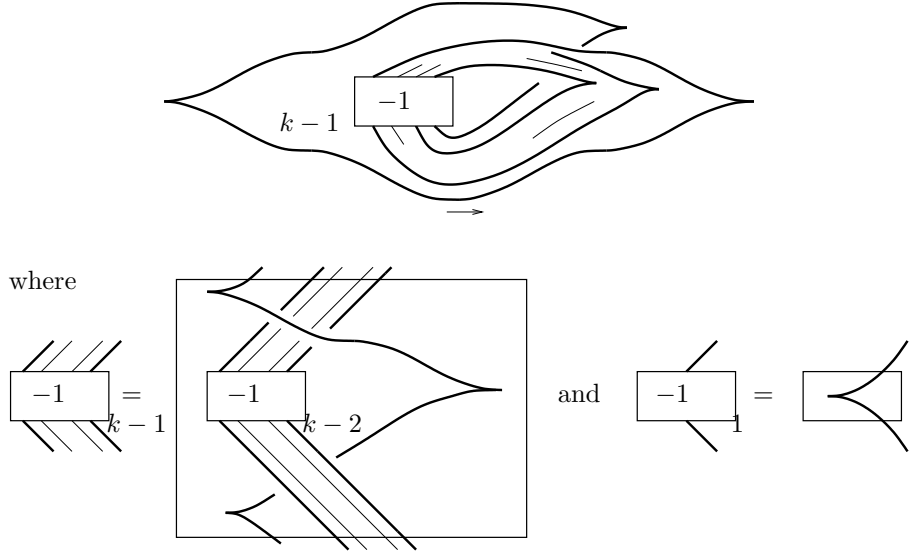}
\extralabel <-3.4cm,2cm> {$-1$}
\extralabel <-11.7cm,2cm> {$-1$}
\extralabel <-9.1cm,2cm> {$-1$}
\extralabel <-7.1cm,5.9cm> {$-1$}
\extralabel <-8.4cm,5.6cm> {$k-1$}
\extralabel <-10.7cm,1.6cm> {$k-1$}
\extralabel <-10.5cm,2cm> {$=$}
\extralabel <-4.7cm,2cm> {and}
\extralabel <-2.1cm,2cm> {$=$}
\extralabel <-2.4cm,1.6cm> {$1$}
\extralabel <-8.1cm,1.6cm> {$k-2$}
\extralabel <-12cm,3.5cm> {where}
\endrelabelbox}
  \caption{The Legendrian knot $K_k$.}
  \label{Kk} 
\end{figure}

In order to have a diagram for each overtwisted contact structure on
$S^1\times S^2$, we first have to find a diagram for contact structures
representing each spin$^c$ structure, and then form the connected sum of
these with the contact structures found in the previous subsection for $S^3$.
Notice that since $H^2(S^1\times S^2; \bfz )\cong \bfz $ has no 2-torsion,
 a spin$^c$ structure is uniquely characterised by its first Chern class.
So the problem reduces to finding a contact structure $\xi _k$ on $S^1\times
S^2$ with $c_1(\xi _k )=2k$ for all $k\in \bfz $. (Recall that the
first Chern class of a 2-plane field is always an even class.)
First we inductively define the Legendrian knot $K_k$ by Figure~\ref{Kk}.

\begin{lemma}
For the oriented Legendrian knot $K_k$ defined by Figure~\ref{Kk},
with $k\geq 2$, we have {\rm {rot}}$(K_k)=k-2$ and ${\rm {tb}}(K_k)=1-k^2$.
\end{lemma}

\begin{proof}
Recall from the proof of Lemma~\ref{lemma:shark}
the formulae for computing tb and rot from the front projection.
Denote the contribution of the box to tb and rot by $t_{k-1}$ and $r_{k-1}$,
respectively. Then by counting the cusps and crossings outside the box we see
\[ \tb (K_k)=t_{k-1}-(k-1)-\frac{1}{2}(k+1)=t_{k-1}-\frac{3}{2}k+\frac{1}{2}\]
and
\[ \rot (K_k)=r_{k-1}+\frac{1}{2}(1-k).\]
From the inductive definition of the box we have the recursive formulae
\[ t_1=-\frac{1}{2},\;\; t_{k-1}=t_{k-2}-2(k-2)-\frac{3}{2}\]
and
\[ r_1=\frac{1}{2},\;\; r_{k-1}=r_{k-2}+\frac{3}{2},\]
from which one finds
\[ t_{k-1}=\frac{1}{2}-k^2+\frac{3}{2}k,\;\;\; r_{k-1}=
-\frac{5}{2}+\frac{3}{2}k.\]
Substituting this into the expressions for tb and rot we obtain the
claimed result.
\end{proof} 

In the following proposition and its proof we use again the notation
of Section~\ref{second}; in particular, $\mu_i$ denotes a meridian
of~$K_i$.

\begin{figure}[h]
\centerline{\relabelbox\small
\epsfxsize 12cm \epsfbox{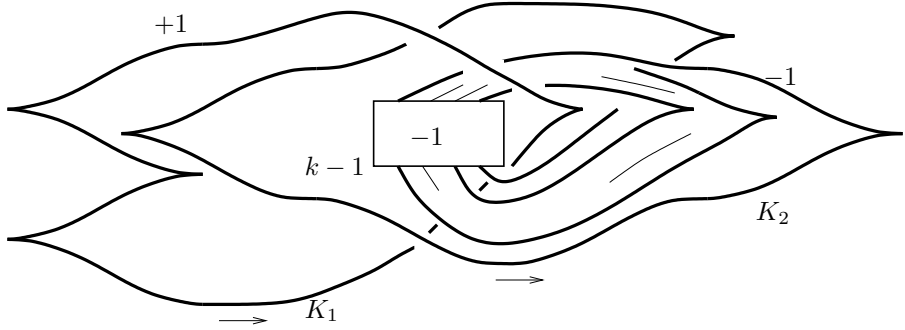}
\extralabel <-6.6cm,2.4cm> {$-1$}
\extralabel <-8.0cm,2.0cm> {$k-1$}
\extralabel <-2.0cm,1.4cm> {$K_2$}
\extralabel <-8.0cm,0.1cm> {$K_1$}
\extralabel <-1.9cm,3.2cm> {$-1$}
\extralabel <-10.0cm,3.9cm> {$+1$}
\endrelabelbox}
  \caption{Surgery diagram of contact structure on $S^1\times S^2$.}
  \label{s1xs2} 
\end{figure}

\begin{proposition}
For $k\geq 2$
the surgery diagram of Figure~\ref{s1xs2} defines a contact structure
$\xi_k$ on $S^1\times S^2$ with $c_1( \xi _k)=(2k-2)PD^{-1}([\mu _2])$.
Here $[\mu_2]$ is a generator of $H_1(S^1\times S^2;{\mathbb Z})\cong
{\mathbb Z}$.
\end{proposition}

Here, by slight abuse of notation, $K_2$ denotes the Legendrian knot $K_k$
considered previously.

\begin{proof}
First of all, we need to check that the topological result of the
described surgeries is $S^1\times S^2$. For that, we observe that
the surgery diagram of Figure~\ref{s1xs2} is topologically
equivalent to that of Figure~\ref{s1xs2a}, where the indicated framings
are now relative to the surface framings in~$S^3$. Blowing down the
$(-1)$-framed unknot $K_1$ (see~\cite[p.~150]{gost99}) adds a $(+1)$-twist
to the $k$ strands of $K_2$ running through it (i.e.\ cancels the
$(-1)$-box) and adds $\lk (K_1,K_2)^2=k^2$ to the framing of $K_2$,
which means that we end up with a single $0$-framed unknot,
which is a surgery picture for $S^1\times S^2$.

\begin{figure}[h]
\centerline{\relabelbox\small
\epsfxsize 8cm \epsfbox{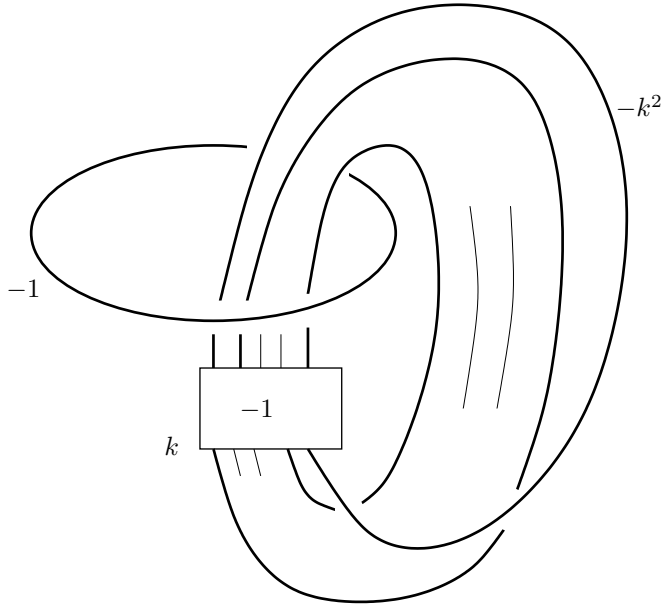}
\extralabel <-5.2cm,2.5cm> {$-1$}
\extralabel <-6.2cm,2cm> {$k$}
\extralabel <-8.3cm,4.1cm> {$-1$}
\extralabel <-.2cm,6.5cm> {$-k^2$}
\endrelabelbox}
  \caption{Surgery diagram for $S^1\times S^2$.}
  \label{s1xs2a} 
\end{figure}

The contact manifold $(S^1\times S^2,\xi_k)$ is the boundary of the almost
complex manifold $(X,J)$ obtained by attaching two $2$-handles to
$D^4$ and forming the connected sum with ${\mathbb C}P^2$ (since
we perform one contact $(+1)$-surgery), in particular, $c_1(\xi_k)$ is the
restriction of $c:=c_1(X,J)$ to the boundary.

Since $\rot (K_1)=1$ and $\rot (K_2)=k-2$, we have (with $[\bfc P^1]$ denoting
the class of a complex line in the $\bfc P^2$ summand)
\[ PD(c)=[N_1]+(k-2)[N_2]+3[\bfc P^1].\]
This implies
\[ c_1(\xi_k)=PD^{-1}([\mu_1]+(k-2)[\mu_2]).\]
With respect to the surface framing in $S^3$, the surgery coefficients
are $n_1=\tb (K_1)+1=-1$ and $n_2=\tb (K_2)-1=-k^2$. Moreover,
we have $\lk (K_1,K_2)=k$. Thus the relations between $[\mu_1]$ and $[\mu_2]$
are given by
\[ -[\mu_1]+k[\mu_2]=0,\;\; k[\mu_1]-k^2[\mu_2]=0.\]
Hence $[\mu_2]$ generates $H_1(S^1\times S^2)$ and
$c_1(\xi_k)=(2k-2)PD^{-1}([\mu_2])$.
\end{proof}

A surgery diagram for an overtwisted contact structure $\xi_0$ on
$S^1\times S^2$ with $c_1(\xi_0)=0$ is given by the disjoint union
of the knots in Figures~\ref{shark} and~\ref{saucer}. (This amounts
to a connected sum of the tight contact structure on $S^1\times S^2$
with an overtwisted contact structure on~$S^3$.)

By rotating the link diagram of Figure~\ref{s1xs2} by $180^{\circ}$
in the plane and keeping the orientations of $K_1$ and $K_2$,
the rotation numbers change sign, while the homology classes
$[\mu_1]$ and $[\mu_2]$ remain unchanged. So this provides
surgery diagrams of contact structures $\xi_{-k}$ on $S^1\times S^2$
with first Chern class $c_1(\xi_{-k})=(2-2k)PD^{-1}[\mu_2]$, $k\geq 2$.

Notice that by reversing the orientations on the knots $K_1$ and
$K_2$ of Figure~\ref{s1xs2} we could achive a sign change in the rotation
numbers, implying a sign change in the coefficient of the expression
for $c_1(\xi _k)$. However, this change would also change the sign of 
$[\mu _2]$, so we would not have gained anything.

\vspace{2mm}

Here, again, is an alternative proof for the construction of all
contact structures on $S^1\times S^2$; we leave it to the reader to check the
details. Let $K_0$ be the Legendrian unknot of Figure~\ref{saucer}
with $\tb (K_0)=-1$. Let $K_1$ be a copy of this knot linked
$k$ times with $K_0$. Let $K_2$ be a Legendrian push-off of $K_1$ with
two zigzags added such that (with the appropriate choice of
orientations) $\rot (K_2)=\rot (K_1)+2=2$. Contact $(+1)$-surgeries on
$K_0,K_1,K_2$ give an overtwisted contact structure on
$S^1\times S^2$ with $c_1=2kPD^{-1}[\mu_0]$, where the class of the
normal circle
$\mu_0$ to $K_0$ generates $H_1(S^1\times S^2)$. That this surgery
picture does indeed, topologically, describe $S^1\times S^2$ can be
seen by sliding $K_2$ over $K_1$.
\subsection*{Overtwisted contact structures on $3$-manifolds}
We now give an algorithm for drawing surgery diagrams for all overtwisted
contact structures on an arbitrary given $3$-manifold~$Y$. Recall from the
discussion at the beginning of this section that we only need to find
diagrams realising all spin$^c$ structures.

Thus, assume that the 3-manifold $Y$ is given by surgery along a
framed link
\[ \Li=((K_1, n_1), \ldots , (K_t, n_t)) \subset S^3. \]
We may assume that these are honest surgeries,
i.e.\ with integer framings~$n_i$. If $Y$ is represented by Dehn surgeries
(with rational coefficients) along a certain link, one can use continued
fraction expansions of the surgery coefficients to turn the diagram into
an integral surgery diagram as above. We retain the notation of
Section~\ref{second}, except that we allow ourselves to identify
the normal circles $\mu_i$ with the homology classes they represent.

\begin{figure}[h]
\centerline{\relabelbox\small
\epsfxsize 8cm \epsfbox{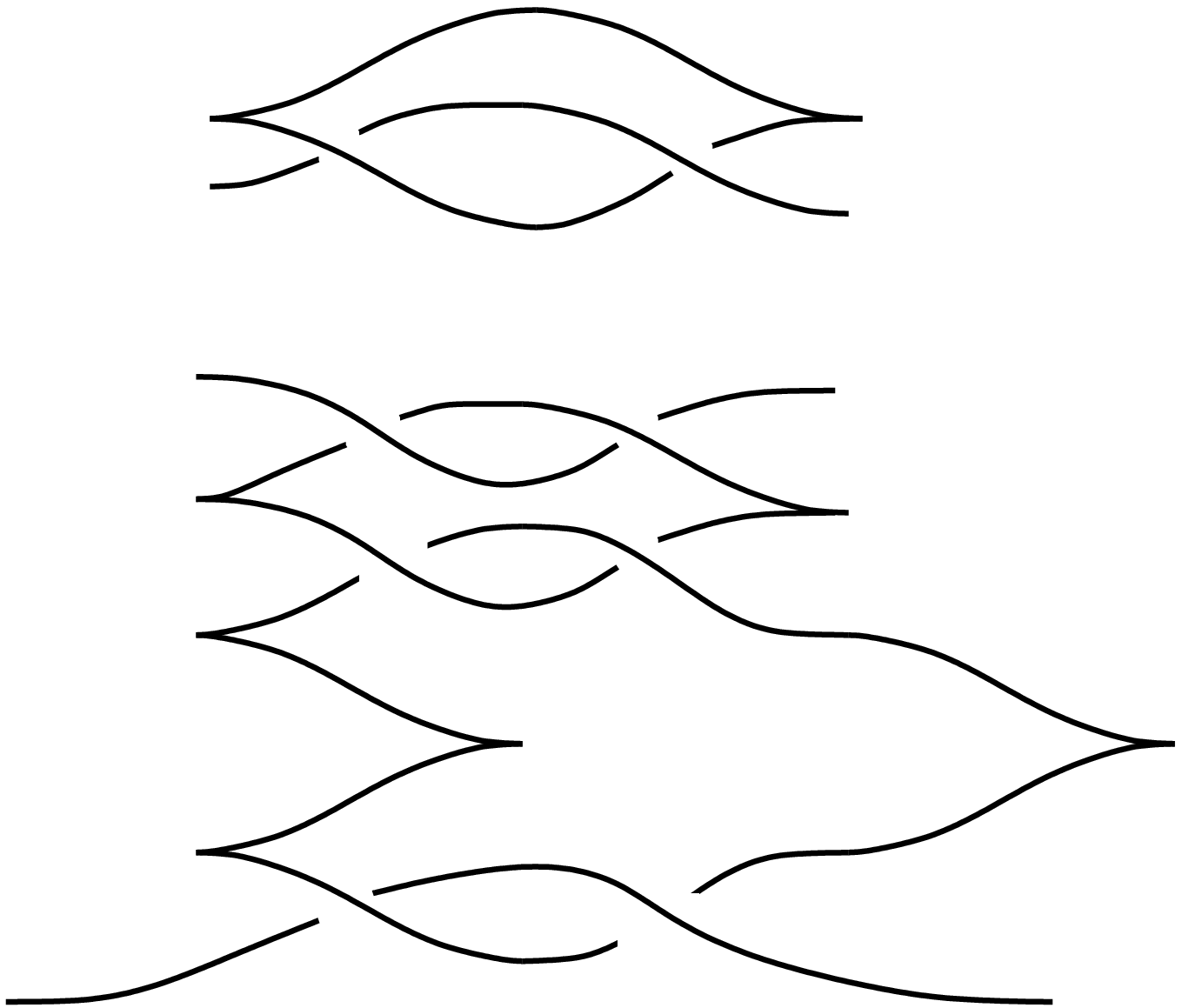}
\extralabel <-4.7cm,4.5cm> {$\vdots$}
\extralabel <-2cm,6cm> {$K_{i,l_i}$}
\extralabel <-7.5cm,6cm> {$-1$}
\extralabel <-2cm,5.3cm> {$K_{i,l_i-1}$}
\extralabel <-7.5cm,5.3cm> {$-1$}
\extralabel <-2cm,3.2cm> {$K_{i,1}$}
\extralabel <-7.5cm,3.2cm> {$-1$}
\extralabel <.2cm,1.7cm> {$K_{i,0}$}
\extralabel <-6cm,1.6cm> {$+1$}
\extralabel <-.4cm,0cm> {$K_i$}
\extralabel <-7.5cm,0.3cm> {$-1$}
\endrelabelbox}
  \caption{The first change on the surgery diagram.}
  \label{change1} 
\end{figure}

In order to find a contact surgery diagram for {\em some} contact
structure on $Y$ we put the knots $K_i$ into Legendrian position
relative to the standard contact structure on $S^3$. Write $b_i$ for
the Thurston-Bennequin invariant $\tb (K_i)$.
If $n_i<b_i$, then by adding zigzags to the Legendrian
knot $K_i$ (which decreases $b_i$) we can arrange
$n_i=b_i-1$, hence contact $(-1)$-surgery on $K_i$ gives the
desired result. If $n_i \geq b_i$ then we transform the Legendrian link
near $K_i$ as shown in Figure~\ref{change1}, where $l_i=n_i-b_i$
and the surgery coefficients have to be read relative to the contact framing. 

Here is the verification that this does indeed correspond to a
surgery along $K_i$ with framing $n_i$ (relative to the surface 
framing in~$S^3$): First of all, we observe that the surgery coefficients
relative to the surface framing are $-2$ for $K_{i,s}$, $s=1,\ldots ,l_i$,
for $K_{i,0}$ it is $-1$, and for $K_i$ it is $b_i-1$.
We now slide off $K_{i,0},K_{i,1},\ldots , K_{i,l_i}$ (in this order).
On sliding off $K_{i,0}$, the topological framing
of $K_{i,1}$ (that is, the framing of the surgery relative to
the surface framing of~$K_{i,1}$)
changes to $-2+1=-1$, that of $K_i$ to $b_i-1+1=b_i$,
and $K_i$ becomes linked once  with $K_{i,1}$. Continuing this way, each step
produces a $(-1)$-framed unknot linked once with $K_i$. Finally, we end
up with $l_i+1$ unknots with topological framing $-1$, which can be
blown down, and with $K_i$ having framing $b_i-1+l_i+1=n_i$, as
desired.

We claim that after the changes described in Figure~\ref{change1}
have been effected, the normal circles to $K_i$, $i=1,\ldots ,t$,
still generate $H_1(\partial X;\bfz )$: Choose orientations on
$K_i, K_{i,0},\ldots ,K_{i,l_i}$ such that the intersection number
of successive knots in this sequence equals $+1$ (this is only necessary to
fix signs in the following computation). Write $\nu_0,\ldots ,\nu_l$
(we suppress the index~$i$) for
the homology classes represented by the normal circles to the knots
$K_{i,0},\ldots ,K_{i,l_i}$. These classes generate $H_1(\partial X;\bfz )$,
and by Section~\ref{second} we have the following relations:
\begin{eqnarray*}
-2\nu_l+\nu_{l-1} & = & 0,\\
-2\nu_i+\nu_{i+1}+\nu_{i-1} & = & 0,\;\; i=1,\ldots ,l-1,\\
-\nu_0+\nu_1+\mu_i & = & 0.
\end{eqnarray*}
The second relation implies $\nu_{i+1}\in\langle\nu_i,\nu_{i-1}\rangle$
for $i=1,\ldots ,l-1$; the third relation yields $\nu_1\in\langle
\nu_0,\mu_i\rangle$. Finally, the relation provided by the Seifert surface
of $K_i$ allows to express $\nu_0$ as a linear combination of
$\mu_1,\ldots ,\mu_t$. In total, we see that all $\nu_{i,j}$ are contained
in the linear span of the $\mu_i$ in $H_1(\partial X;\bfz )$.

We have thus found a contact $(\pm 1)$-surgery description for
some contact structure on the given manifold~$Y$. We now should like
to perform further changes on that surgery diagram so as to
realise all possible spin$^c$ structures. The idea behind the
following construction is first to introduce additional
surgery curves such that (a)~appropriate surgeries along these
curves do not change the topology of $Y$ and (b)~a subset of
the additional surgery curves corresponds to a description of $S^1\times S^2$.
Then the ideas used previously for $S^1\times S^2$ can be applied again.

\begin{figure}[h]
\centerline{\relabelbox\small
\epsfxsize 10cm \epsfbox{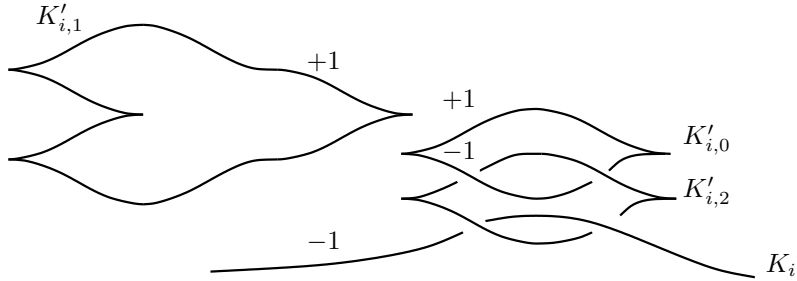}
\extralabel <-9.6cm,3.4cm> {$K_{i,1}'$}
\extralabel <-1cm,1.8cm> {$K_{i,0}'$}
\extralabel <-1cm,1.1cm> {$K_{i,2}'$}
\extralabel <.1cm,0.1cm> {$K_i$}
\extralabel <-6cm,0.4cm> {$-1$}
\extralabel <-4.2cm,1.6cm> {$-1$}
\extralabel <-4.2cm,2.3cm> {$+1$}
\extralabel <-6cm,2.8cm> {$+1$}
\endrelabelbox}
  \caption{The reference contact structure on~$Y$.}
  \label{change2} 
\end{figure}

Consider the contact manifold obtained by adding, for each $i=1,\ldots ,t$,
three surgery curves
$K_{i,0}', K_{i,1}', K_{i,2}'$ as
indicated in Figure~\ref{change2}.

Observe that the topological framings of 
$K_{i,0}', K_{i,1}', K_{i,2}'$ are $0$, $-1$, and $-2$, respectively.
Hence, with appropriate orientations on these knots and with
$\mu_0',\mu_1',\mu_2'$ denoting the homology classes represented by
the normal circles to these knots (again we suppress the index~$i$), we
have the relations
\begin{eqnarray*}
0\cdot\mu_0'+\mu_2' & = & 0,\\
-\mu_1' & = & 0,\\
-2\mu_2'+\mu_0'-\mu & = & 0,
\end{eqnarray*}
that is, $\mu_{i,0}'=\mu_i$ and $\mu_{i,1}'=0=\mu_{i,2}'$, $i=1,\ldots ,t$.
Observe that the surgery curve $K_{i,0}'$ on its own gives a description
of $S^1\times S^2$, with first homology group generated by $\mu_{i,0}'$.

Topologically, these additional surgery curves do not change anything,
so that we still have a description of~$Y$: The $(-1)$-framed unknot
$K_{i,1}'$ gives a trivial surgery; a slam-dunk of $K_{i,0}'$ changes
the framing of $K_{i,2}'$ to $\infty$, which again gives a trivial
surgery. The presence of $K_{i,1}'$ ensures that the diagram describes
an {\em overtwisted} contact structure $\xi_0$ on~$Y$, which will be
our reference contact structure, inducing the spin$^c$ structure
$\t_0=\t_{\xi_0}$.

By viewing the knots in this diagram as attaching circles of $2$-handles
rather than surgery curves, we can read the diagram as a description
of a $4$-manifold $X$ with boundary~$Y$. We have seen that, away from
finitely many points, $X$ admits an almost complex structure $J$ such that
$\xi_0=T\partial X\cap J(T\partial X)$. The corresponding spin$^c$
structure $\s_0$ on $X$ restricts to $\t_0$ along $Y=\partial X$.

Given $\t\in Spin^c(Y)$ there is, thanks to the free and transitive action
of $H^2(Y;\bfz )$ on $Spin^c(Y)$, a class $a_{\t}\in H^2(Y;\bfz )$
such that $\t =\t_0\otimes a_{\t}$. Since the restriction homomorphism
$H^2(X;\bfz )\rightarrow H^2(Y;\bfz )$ is surjective (under Poincar\'e
duality this corresponds to the surjectivity of $\varphi_2$ in
Section~\ref{second}), we may assume that $a_{\t}$ lives in
$H^2(X;\bfz )$. Then $\s_0\otimes a_{\t}$ is a spin$^c$ structure
on $X$ that on $Y$ restricts to~$\t$. The advantage of working
over $X$ is that due to $\pi_1(X)=0$ the first Chern class captures
the spin$^c$ structure, whereas on $Y$ the identification of
spin$^c$ structures is complicated by the possible presence of $2$-torsion.

\begin{figure}[h]
\centerline{\relabelbox\small
\epsfxsize 10cm \epsfbox{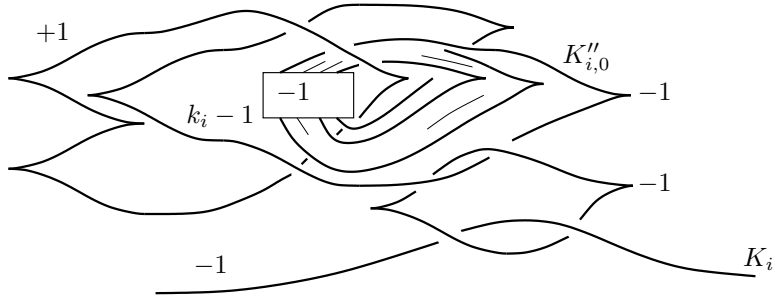}
\extralabel <-9.6cm,3.4cm> {$+1$}
\extralabel <-6.4cm,2.6cm> {$-1$}
\extralabel <-1.6cm,2.6cm> {$-1$}
\extralabel <-2.6cm,3.1cm> {$K_{i,0}''$}
\extralabel <-1.6cm,1.4cm> {$-1$}
\extralabel <-.2cm,.4cm> {$K_i$}
\extralabel <-7.5cm,.3cm> {$-1$}
\extralabel <-7.6cm,2.3cm> {$k_i-1$}
\endrelabelbox}
  \caption{The surgery diagram for $k_i=\alpha_i+1\geq 2$.}
  \label{change3} 
\end{figure}

In conclusion, we need to find a contact surgery diagram that topologically
yields $Y$ and such that the induced spin$^c$ structure $\s$ on $X$
satisfies $c_1(\s )=c_1 (\s_0)+2a_{\t}\in H^2(X;\bfz )$. Observe that
because of $\mu_{i,0}'=\mu_i$, we can --- with $N_{i,0}'$ denoting the normal
disc bounded by $\mu_{i,0}'$ --- write $a_{\t}$ as
\[ a_{\t}=\sum_{i=1}^t\alpha_iPD^{-1}[N_{i,0}']\in H^2(X;\bfz ).\]
If $\alpha_i=0$, we retain the diagram of Figure~\ref{change2} near~$K_i$.
If $\alpha_i>0$, we use instead the diagram depicted in Figure~\ref{change3},
which is modelled on the one we used for $S^1\times S^2$.

Observe that the presence of the (contact) $(+1)$-framed unknot with
Thurston-Bennequin invariant $-2$ (and the fact that the other link
components may be assumed not to intersect the overtwisted disc we
exhibited in Figure~\ref{ot}) again ensures that the resulting
contact structure is overtwisted. Moreover, the diagram is topologically
equivalent to the one of Figure~\ref{change2}, with $K_{i,0}''$
taking the role of~$K_{i,0}'$. Thus, a calculation completely analogous
to the one above for the contact structure $\xi_k$ on
$S^1\times S^2$ shows that passing from the diagram in Figure~\ref{change2}
to the one in Figure~\ref{change3} adds a summand $(2k_i-2)PD^{-1}[N_{i,0}']
=2\alpha_iPD^{-1}[N_{i,0}']$ to the first Chern class of the
corresponding spin$^c$ structure.

For $\alpha_i<0$ one argues similarly, using the diagrams for the
$\xi_{-k}$ instead. This concludes the construction of surgery diagrams
for all overtwisted contact structures on the given~$Y$.

\vspace{1mm}

Notice that when we claim to have found surgery diagrams for {\em all}
overtwisted contact structures on a given (closed) $3$-manifold $Y$, we do of
course rely on Eliashberg's result~\cite{elia89} that overtwisted
contact structures which are homotopic as $2$-plane fields are in fact
isotopic as contact structures. However, our argument clearly provides
an independent proof of the Lutz-Martinet theorem:

\begin{corollary}[Lutz-Martinet]
On any given closed, orientable $3$-manifold, each homotopy class
of $2$-plane fields contains an (overtwisted) contact structure.
\hfill $\Box$
\end{corollary}

For an exposition of the original proof of that theorem, based on surgery
along curves {\em transverse} to a given contact structure, see~\cite{geig}.

\section{$(+1)$-surgery revisited}
\label{section:+1}
In this final section we briefly return to the issues raised in
Remarks \ref{rem:pm} and \ref{rem:tb} concerning the extension of the almost
complex structure over the handle and the value of $c(\Sigma )$ in
the case of contact $(+1)$-surgery. In fact, most of our discussion in the
present section relates to the translation from Weinstein's description
of contact surgery via symplectic handlebodies with contact type boundary
to Eliashberg's description via Stein manifolds (or complex handlebodies with
strictly pseudoconvex boundary), and thus it applies equally well to
the case of contact $(-1)$-surgery. Specifically, we address the question
how to deform a handle in Weinstein's picture so that the contact
structure on the boundary of the handle is given by almost complex
tangencies; we are not concerned with the more subtle point of
the integrability of that almost complex structure (extending
a given complex structure on the initial handlebody). The second issue then
is to give a geometric description for the obstruction to extending
that almost complex structure over the full handle in the case
of contact $(+1)$-surgery -- in the case of $(-1)$-surgery there is
no such obstruction, as already discussed. We hope that the
following considerations will prove useful in other instances where it may
be opportune to switch between Eliashberg's and Weinstein's description
of contact surgery.

We begin with the following simple lemma:

\begin{lemma}
\label{lem:ac}
Let $E\rightarrow X$ be an oriented $\bfr^4$-bundle (over some
manifold~$X$) with bundle metric $g$ and $\xi\subset E$ an oriented
$\bfr^2$-subbundle. Then there is a unique complex bundle structure $J$
on $E$ such that
\begin{itemize}
\item[{\it (i)}] $g$ is $J$-invariant.
\item[{\it (ii)}] $\xi$ is $J$-invariant
\item[{\it (iii)}] $J$ induces the given orientations of $E$ and~$\xi$.
\end{itemize}
Any two complex bundle structures $J_0, J_1$ on $E$ satisfying
(ii) and (iii) are homotopic.
\end{lemma}

\begin{proof}
Let $(\sigma_1,\sigma_2,\sigma_3,\sigma_4)$ be an ordered quadruple
of local $g$-orthonormal sections of $E$ with $(\sigma_1,\sigma_2)$
sections of $\xi$, inducing the given orientations. Then $J$ with the described
properties can be defined by $J\sigma_1=\sigma_2$ and $J\sigma_3=\sigma_4$,
and it is a straightforward check that this is the only way to define~$J$.

Given $J_0,J_1$ as described, let $g_i$, $i=0,1$, be a $J_i$-invariant
bundle metric on~$E$. The first part of the proof tells us that
$J_i$ can be recovered from~$g_i$. The complex bundle structure
$J_t$ corresponding in this way to the bundle metric $g_t=(1-t)g_0
+tg_1$, $t\in [0,1]$, defines a homotopy between $J_0$ and~$J_1$.
\end{proof}

\begin{figure}[h]
\centerline{\relabelbox\small
\epsfxsize 8cm \epsfbox{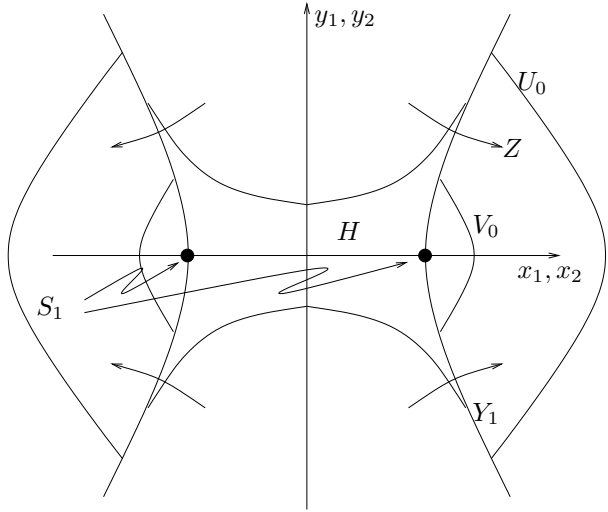}
\extralabel <-7.6cm,2.6cm> {$S_1$}
\extralabel <-1.4cm,4.7cm> {$Z$}
\extralabel <-1.2cm,3.1cm> {$x_1,x_2$}
\extralabel <-3.9cm,6.5cm> {$y_1,y_2$}
\extralabel <-1.8cm,1.2cm> {$Y_1$}
\extralabel <-3.6cm,3.6cm> {$H$}
\extralabel <-1.2cm,5.6cm> {$U_0$}
\extralabel <-1.8cm,3.7cm> {$V_0$}
\endrelabelbox}
  \caption{Contact $(+1)$-surgery}
  \label{figure:+1surgery} 
\end{figure}

Recall from \cite[Section~3]{dige} the description of contact $(+1)$-surgery
as a symplectic handlebody surgery on the concave end of a symplectic
cobordism: Consider $\bfr^4$ with cartesian coordinates $(x_1,y_1,x_2,y_2)$
and standard symplectic form
\[ \omega =dx_1\wedge dy_1+dx_2\wedge dy_2.\]
Then
\[ Z=2x_1\partial_{x_1}-y_1\partial_{y_1}+2x_2\partial_{x_2}-
y_2\partial_{y_2}\]
is a Liouville vector field for $\omega$, that is, ${\mathcal L}_Z\omega
=\omega$. This implies that $\alpha=i_Z\omega$ is a contact form on
any hypersurface transverse to~$Z$. Let $f\colon\bfr^4\rightarrow\bfr$ be the
function defined by
\[ f(x_1,y_1,x_2,y_2)=x_1^2-\frac{1}{2}y_1^2+x_2^2-\frac{1}{2}y_2^2\]
and set $Y_{\mu}=\{ f=\mu\}$ and $S_1=
Y_1\cap\{ y_1=y_2=0\}$, which is Legendrian in $(Y_1,\ker\alpha )$.
A neighbourhood of $S_1$ in $Y_1$ can be
identified with a neighbourhood of a given Legendrian knot $K$ in
$(S^3,\xi_{st})$ (which we take to be the boundary of $D^4$ with its
standard complex structure~$J$),
and Figure~\ref{figure:+1surgery}
shows how to attach a symplectic handle $H$ along $S_1\equiv K$.
(More generally, one can assume that $K$ is a Legendrian knot in
a contact manifold $(Y,\xi )$ given as the boundary of an almost
complex manifold $(X,J)$.)
In \cite{dige} the framing of this surgery is computed to be indeed
$+1$ with respect to the contact framing of~$K$.

The orientation of $Y_1$ is given by $\alpha\wedge d\alpha =i_Z\omega^2/2$.
Hence, in order for $Y_1$ to carry the boundary orientation of
$X_1=\{ f\geq 1\}$, we need to equip $\bfr^4$ with the orientation
given by $-\omega^2$ (or $-df\wedge\alpha\wedge d\alpha$).

A complex bundle structure $J_0$ on $E=T(\bfr^4-\{ 0\} )$ is defined,
in the sense of the preceding lemma, by the $2$-plane bundle
\[ \xi_0 =\ker df\cap\ker\alpha\]
(oriented by $d\alpha$) and the standard metric $g_0$ on $\bfr^4$. Then on
each level surface $Y_{\mu}$ (except at the singular point $0\in Y_0$),
the contact structure $\xi_0$ coincides with the $J_0$-complex
tangencies of~$Y_{\mu}$.

\begin{proposition}
In the notation of Proposition~\ref{prop:+1}, we have $c(\Sigma )=\rot (K)$,
independently of the value of $\tb (K)$.
\end{proposition}

\begin{proof}
We should like to argue that $J_0$ does in fact define
the extension $J$ of the almost complex structure on $D^4$ over $H-\{ 0\}$.
Unfortunately,
this is not quite true, since the boundary of $H$ is not a level surface
of $f$, so the contact structure $\ker\alpha\cap T(\partial H)$ on
$\partial H$ does not
coincide with $\xi_0$, i.e.\ that contact structure is not given by
the $J_0$-complex tangencies of $\partial H$.
Up to homotopy, however,
this is essentially true. Thus, before addressing this mild subtlety,
we prove that $c(\Sigma )=\rot (K)$ from the $J=J_0$ as described.

Let $F$ be the Seifert surface of $K$ in $S^3$ and $D$ the core disc
of $H$,
\[ D =\{ (x_1,y_1,x_2,y_2)\colon
x_1^2+x_2^2\leq 1,\; y_1=y_2=0\} , \]
perturbed slightly around $0$ so that it stays inside $H$ but
misses the origin of $\bfr^4$. Then $\Sigma$, by definition, is
the surface obtained by gluing $F$ and $D$ along $S_1\equiv K$, with
orientation of $K$ equal to the boundary orientation of~$F$.

Along $S^3$ the tangent bundle of $D^4$
splits (as a complex bundle) into the complex line bundle $\xi_{st}$
and a trivial complex line bundle defined by the complex
lines containing the outward normal.
That latter trivialisation extends to a trivialisation of a
complex line bundle in $T\bfr^4|_D$ complementary to $\xi_0 |_D$,
viz., the $J_0$-complex lines containg $Z$.
Therefore the first Chern class $c$ of $J$, when restricted to $\Sigma$, equals
the first Chern class of $\xi|_{\Sigma}$ (with $\xi=\xi_{st}$ on
$F$ and $\xi=\xi_0$ on~$D$).

Moreover, the vector field
\[ v= 2x_2\partial_{x_1}+y_2\partial_{y_1}-2x_1\partial_{x_2}
-y_1\partial_{y_2} \]
is a nowhere zero vector field in $\xi_0|_{H-\{ 0\} }$ -- in particular,
it defines a trivialisation
of the complex line bundle $\xi_0 |_D$ --
and its restriction to $S_1$ is tangent to that circle. By our orientation
assumption on $K\equiv S_1$ and $F$, the value $c(\Sigma )=
\langle c_1(\xi|_{\Sigma}),[\Sigma ]\rangle$ is equal to
the rotation number of $v|_K$ relative to a trivialisation of $\xi_{st}|_F$,
which by definition is precisely $\rot (K)$.
\end{proof}

We now show how to deform the local picture of Figure~\ref{figure:+1surgery}
in such a way that the extension of the almost complex structure over
$H-\{ 0\}$ is indeed defined by~$J_0$.

First of all, we have a contactomorphism $\varphi$ from a neighbourhood
of $K$ in $(S^3,\xi_{st})$ to a neighbourhood of $S_1$ in $(Y_1,\xi_0)$.
Extend $\varphi$ to a diffeomorphism of a neighbourhood of $K$ in $D^4$
to a neighbourhood of $S_1$ in~$X_1$. We claim
that one can homotope $J$ on $D^4$
to an almost complex structure (still denoted~$J$) such that
\begin{itemize}
\item $\xi_{st}$ is still given by the $J$-complex tangencies of $S^3=
\partial D^4$,
\item the homotopy is supported in a given neighbourhood of $K$ in~$D^4$,
\item $\varphi_*J$ coincides with $J_0$ in a neighbourhood $U_0$ of
$S_1$ in $X_1$.
\end{itemize}

In order to see this, extend $\xi_{st}$ to a plane field $\xi$ on
$D^4-\{ 0\}\subset\bfc^2$ as the complex tangencies of the spheres of
radius $r\in (0,1]$. Since $\xi_{st}$ coincides with $\varphi^*\xi_0$
on a neighbourhood of $K$ in~$S^3$, there is a homotopy of~$\xi$, fixed
on $S^3$ and supported in a neighbourhood of $K$ in~$D^4$, to a plane
field (still denoted~$\xi$) that coincides with $\varphi^*\xi_0$ in
a (smaller) neighbourhood $U$ of $K$ in~$D^4$. Clearly, there is a
corresponding homotopy of the standard metric on $D^4$ to a metric coinciding
with $\varphi^*g_0$ near~$K$. Lemma~\ref{lem:ac} then allows us to
construct the desired homotopy of~$J$, with $U_0=\varphi (U)$.

We attach the handle $H$ inside the neighbourhood $U_0\cap Y_1$. Next choose a
smaller neighbourhood $V_0\subset U_0$
of $S_1$ in $X_1$ such that $V_0\cap Y_1$
lies completely inside the region where $H$ is attached to~$X_1$.
Let $h_0\colon \bfr^4-\{ 0\}\rightarrow\bfr^-$ be the function
\[ h_0(x_1,y_1,x_2,y_2)=-\frac{2}{4x_1^2+y_1^2+4x_2^2+y_2^2},\]
and $h\colon\bfr^4\rightarrow\bfr^-$ a smooth function such that
$h=h_0$ outside a neighbourhood of the origin chosen so small that
the flow $\varphi_t$ of $hZ$ coincides with the flow of
$h_0Z$ on a collar neighbourhood $W_0$ of $Y_1-\overline{V_0}$ in~$X_1$.
Notice that since the flow $\varphi_t$ of $hZ$ is simply a reparametrisation of
the flow of~$Z$, hypersurfaces transverse to $Z$ stay transverse to $Z$ and
continue to inherit a contact structure from the $1$-form $\alpha =i_Z\omega$.

Observe that ${\mathcal L}_{hZ}\alpha=i_{hZ}d\alpha=h\alpha$, so the flow
$\varphi_t$ of $hZ$ preserves $\ker\alpha$. Furthermore, $df(h_0Z)\equiv -2$.
This implies that $\varphi_t(Y_1-\overline{V_0})\subset Y_{1-2t}$ and
\[ \varphi_{t*}(\ker df_x\cap\ker\alpha_x)=\ker df_{\varphi_t(x)}\cap
\ker\alpha_{\varphi_t(x)}\;\;\mbox{\rm for}\; x\in W_0,\]
in particular, the map $\varphi_t\colon Y_1-\overline{V_0}
\rightarrow Y_{1-2t}$ is an embedding
preserving the contact structure $\xi_0$ on the respective hypersurfaces.

So $\varphi^*\varphi_t^*\xi_0$ (on~$\varphi^{-1}(U_0)=U$) is a homotopy
of $\varphi^*\xi_0=\xi$ that stays
constant in the collar neighbourhood $\varphi^{-1}(W_0)$ of $S^3\cap
\varphi^{-1}(U_0-\overline{V_0})$. This allows to spread out that
homotopy over a collar of $S^3$ in $D^4$ so as to obtain a plane field
(still denoted~$\xi$) on $D^4-\{ 0\}$ that is homotopic to the old
$\xi$ under a homotopy supported in a neighbourhood of $K$ in~$D^4$.
Once again, Lemma~\ref{lem:ac} defines a corresponding homotopy of~$J$
(since one can always interpolate between different metrics).

Thus, after such a homotopy of $J$
and a homotopy of $\xi_{st}$ defined by $\ker\varphi^*\varphi_t^*
\alpha |_{TS^3}$, fixed outside $S^3\cap\varphi^{-1}
(\overline{V_0})$,
we may assume that $\varphi_1\circ\varphi$ sends
$S^3\cap\varphi^{-1}(U_0)$ contactomorphically into
$\varphi_1(Y_1)$ and that $\varphi_1\circ\varphi$ is a
$J$-$J_0$-holomorphic map on a collar neighbourhood of $S^3\cap\varphi^{-1}
(U_0)$ in~$D^4$. Notice, however, that $\ker\varphi^*\varphi_1^*
\alpha |_{TS^3}$
need no longer coincide on $\varphi^{-1}(\overline{V_0})$ with the (homotoped)
$J$-complex tangencies, and $(\varphi_1
\circ\varphi )_*\xi_{st}$ may not coincide with the $J_0$-complex tangencies
of $\varphi_1(Y_1\cap \overline{V_0})$.

Now define $H'$ to be the region bounded by $Y_{-1}$ and $\varphi_1(Y_1)$;
this really amounts to a deformation of $\varphi_1(H)$ keeping its boundary
transverse to~$Z$, hence to a contact isotopy of the surgered contact manifold.
This $H'$ defines contact $(+1)$-surgery in such a way that the extension
of the almost complex structure over $H'-\{ 0\}$ is defined by~$J_0$.

\vspace{2mm}

Finally, we want to give a more geometric argument for the extendability
of the almost complex structure $J_0$ on $H-\{ 0\}$ to $H\#\bfc P^2$;
this gives a new proof of the statement $d_3(\xi_H)=1/2$ in
Proposition~\ref{prop:+1}, independently of the value of $\tb (K)$.

To that end, consider the map $\pi\colon\bfr^4\rightarrow\bfc$ given by
$\pi (x_1,y_1,x_2,y_2)=z_1^2+z_2^2$, where we set
$z_1=x_1+iy_1$ and $z_2=x_2+iy_2$. Write $\pi_1,\pi_2$ for the real and
imaginary part of $\pi$, respectively, i.e.
\begin{eqnarray*}
\pi_1(x_1,y_1,x_2,y_2) & = & x_1^2-y_1^2+x_2^2-y_2^2,\\
\pi_2(x_1,y_1,x_2,y_2) & = & 2x_1y_1+2x_2y_2.
\end{eqnarray*}
Then
\begin{eqnarray*}
d\pi_1 & = & 2x_1dx_1-2y_1dy_1+2x_2dx_2-2y_2dy_2,\\
d\pi_2 & = & 2x_1dy_1+2y_1dx_1+2x_2dy_2+2y_2dx_2.
\end{eqnarray*}

There is an obvious linear homotopy on $\bfr^4-\{ 0\}$ between the pair
$(df,\alpha )$ and the pair $(d\pi_1,d\pi_2)$, the homotopy being through
linearly independent pairs of $1$-forms. Therefore, $J_0$ is homotopic,
by Lemma~\ref{lem:ac}, to the almost complex structure $J_1$ determined by
the plane field $\ker d\pi_1\cap\ker d\pi_2$, coorientation given by
$-d\pi_1\wedge d\pi_2$, and ambient orientation given by $-\omega^2$.
This $J_1$ is exactly the almost complex structure near an
incorrectly oriented critical point (excluding that point) of an achiral
Lefschetz fibration, see~\cite[Section~8]{gost99}, and Lemma~8.4.12
of the cited reference
provides a geometric argument, based on work of Matsumoto, for the
extendability of $J_1$ over the connected sum with a copy of $\bfc P^2$.

\vspace{2mm}

\noindent {\bf Acknowledgements.}
F.~D.\ is partially supported by grant no.\ 10201003 of the National Natural
Science Foundation of China. H.~G.\ is partially supported by grant no.\
GE 1254/1-1 of the Deutsche Forschungsgemeinschaft within the
framework of the Schwerpunktprogramm 1154 ``Globale
Differential\-geo\-metrie''. A.~S.\ is partially supported by
OTKA T034885.

H.~G.\ and A.~S.\ acknowledge the support of TUBITAK for attending
the 10th G\"okova Geometry-Topology conference, which allowed us to discuss
some of the final details of this paper.

\end{document}